\documentclass{article}
\usepackage{amsmath,amscd,amssymb,amsthm}
\input xy
\xyoption{all}

\newtheorem{thm}{Theorem}[section]
\newtheorem{prop}[thm]{Proposition}
\newtheorem{lem}[thm]{Lemma}

\newtheorem{cor}[thm]{Corollary}

\theoremstyle{definition}

\newtheorem{rem}{Remark}

\begin{document}
\setlength{\parskip}{.25cm}
\title{The Singularities of the Wave Trace of the Basic Laplacian of a Riemannian Foliation}
\author{M. R. Sandoval\\ Department of Mathematics\\ 
Trinity College\\ Hartford, Connecticut 06106\\
\\
mary.sandoval@trincoll.edu\\
(telephone) 860-297-2016\\
(fax) 860-987-6239}

\maketitle

\pagebreak

\begin{abstract}
We apply techniques of microlocal analysis to the study of the transverse geometry of Riemannian foliations in order to analyze spectral invariants of the basic Laplacian acting on functions on a Riemannian foliation with a bundle-like metric.  In particular, we consider the trace of the basic wave operator when the mean curvature form is basic. We extend the concept of basic functions to distributions and demonstrate the existence of the basic wave kernel.  The singularities of the trace of this basic wave kernel occur at the lengths of certain geodesic arcs which are orthogonal to the closures of the leaves of the foliation. In cases when the foliation has regular closure, a complete representation of the trace of the basic wave kernel can be computed for $t\not=0$.  Otherwise, a partial trace formula over a certain set of lengths of  well-behaved geodesic arcs is obtained.
\end{abstract}


\pagebreak

%
%
%

\renewcommand{\baselinestretch}{2}

\begin{section}{Introduction}
Let $M$ be a compact manifold without boundary of dimension $n$ that admits a Riemannian foliation $\mathcal{F}$ of  dimension $p.$ Let $q$ denote the codimension of the foliation. Recall that a foliation is Riemannian if there is a metric on $M$ with respect to which the distance between leaves is locally constant.  Such a metric is said to be bundle-like with respect to the foliation $\mathcal{F}.$ We will assume that $(M,\mathcal{F})$ is equipped with just such a metric. We will denote the individual leaves of this foliation by $L$, and the associated distribution by $T\mathcal{F}\subset TM.$  Note that we have the following short exact sequence:
\begin{equation}
0\rightarrow T\mathcal{F}\rightarrow TM\rightarrow Q:=TM/T\mathcal{F}\rightarrow 0.
\end{equation}
If the metric is bundle-like, we have an isomorphism $Q\cong (T\mathcal{F})^\perp=N\mathcal{F},$ and the metric induces a transverse metric $g_T$ on $Q$.  Conversely, given a transverse metric $g_T,$ there exist bundle-like metrics on $M$ which have $g_T$ as their associated transverse metric.  In this paper, we are interested in the manner in which the leaves are glued together to form the manifold $M$.   This is (roughly) the transverse geometry of the foliation.  Broadly speaking, the theme of this paper is to determine to what extent one can associate geometric objects on $Q^*$ with analytic objects that are associated to the transverse structure of the foliation.

An important class of functions that are associated to the transverse structure of $(M,\mathcal{F})$ are the basic functions; these are the functions on $M$ that are constant along the leaves of the foliation, denoted by $C^\infty_B(M,\mathcal{F}).$ (Note: if the foliation contains a dense leaf, or a leaf that is always contained in the closure of any other leaf (like the Reeb foliation) then the basic functions are just the constant functions.  We will focus on the opposite case--the case where the set of basic functions is infinite dimensional.)  Observe that if a function is basic, it is also constant on the closures of the leaves. (In general, the leaves themselves may not be closed.)  In fact, the dimension of the closures of the leaves of an arbitrary foliation may vary over $M.$ Thus, the partition of $M$ into leaf closures may not form another foliation of $M$.  It does, however, have a nice structure--that of a {\it singular Riemannian foliation}.  (See Chapter 6 of \cite{Molino} for definitions.)

There is a similar notion of being basic that applies to forms: a form $\alpha\in\Omega^k(M)$ is said to be basic if $i_X\alpha=i_Xd\alpha=0$ for every $X\in C^\infty(T\mathcal{F}).$ Of particular interest for the purposes of this paper is the mean curvature 1-form, $\kappa,$ given by:
\begin{equation}
\kappa(Z)=\sum_{i=1}^pg(\nabla_{E_i}E_i, Z),\quad \text{where }Z\in C^\infty(N\mathcal{F}) 
\end{equation}
(see \cite{T1}).  (In the above, $g$ denotes the metric on $M$, and the $E_i$, $i=1,\dots p$ are a basis of $T\mathcal{F}.$) This notion turns out to be important in defining a version of the Laplacian on basic functions.

The ordinary Laplacian, $\Delta,$ with respect to an arbitrary bundle-like metric $g$ does not, as a general rule, preserve the space of basic functions.  However, one can define an associated operator on the space of basic functions (and also, incidentally, on the space of basic forms), called the basic Laplacian.  The basic Laplacian, $\Delta_B,$ is equal to $\delta_Bd_B+d_B\delta_B$ where $d_B$ is the exterior derivative restricted to basic functions (or forms) and $\delta_B$ is its adjoint.  It has been shown in \cite{PaRi} that the ordinary Laplacian $\Delta$ restricts to $\Delta_B$ precisely when $\kappa$ is a basic 1-form.  Thus, if $\kappa$ is basic, then the spectrum of $\Delta_B$ is contained in the spectrum of the ordinary Laplacian, $\Delta.$  In fact, the authors of \cite{PaRi} have shown that there is a natural projection $P$ from $C^\infty(M)\rightarrow C^\infty_B(M,\mathcal{F}),$ the basic projection, and that 
\begin{equation}\label{e:laps}
\Delta_BP=P\Delta.
\end{equation}

It is always possible to find a bundle-like metric for which $\kappa$ is basic by the results of \cite{Dom}. It is even possible to pick a bundle-like metric for which $\kappa$ is basic and the induced transverse metric $g_T$ is prescribed, \cite{LR2}. However, the spectrum of $\Delta_B$ depends on both transverse and leaf-wise properties of the given bundle-like metric.  In particular, it has been shown in \cite{KR2} that the eigenvalues of basic Laplacian depend on the volumes of the leaf closures, and thus, the basic spectrum depends on the choice of the entire bundle-like metric and not just the transverse part.

The goal of this paper is to compute invariants of the basic spectrum in terms of the global structure of the foliation.  The approach we will use in this paper is via the kernel to the wave operator for the basic Laplacian, in the spirit of \cite{DG}.  Many of the results follow from straightforward application of the results of \cite{H3}, \cite{DG}, and \cite{Z}. Recall that the wave equation admits a fundamental solution in the category of generalized functions. Thus, we first show that the notion of basic functions can be extended to distributions (in the analytic sense). (Note: in what follows, we will refer to distributions in the analytic sense as generalized functions, and reserve the term distribution for the association of vector subspaces of $T_pM$ to points in $p\in M$ in cases where confusion may result.)  We then define the basic wave kernel in an analogous manner to the basic heat kernel, and demonstrate that the basic wave kernel exists and is related to the ordinary wave kernel via the basic projection $P$.  We show that the singularities of the trace of the basic wave kernel are contained in the set of lengths of certain geodesics arcs which are orthogonal to the leaf closures. These lengths are invariants of the basic spectrum.  Furthermore, if one can localize to avoid certain particularly problematic values of $T$, one can derive a representation of the basic wave trace as a sum of Lagrangian generalized functions on $T^*\mathbb{R}.$ If the foliation admits regular closure (that is, when the closures of the leaves of the foliation all have the same dimension),  the trace of the basic wave kernel has a representation as a sum of Lagrangian generalized functions near any singularity.  

The heat kernel and the basic spectrum have been widely studied by many researchers, including \cite{KR2}, \cite{NRT}, \cite{LR1}, \cite{LR2}, \cite{Ri2}.  In particular, in \cite{KR2}, the researcher showed that the basic heat kernel pulled back to the diagonal in $M\times M:$  $K_B(t,x,x)$ admits an asymptotic expansion in $t$ whose coefficients are invariants of the basic spectrum.  These coefficients depend on the codimension of the leaf closures, the volumes of the leaf closures in $M$ and the lifted foliation on the oriented orthornormal transverse frame bundle, $\widehat{M}$, the curvature at $x\in M$ and the curvature of a related manifold, the basic manifold, $W$. (See Theorem 3.1 of \cite{KR2} for precise statements.)  One feature of the asymptotic formula for the heat kernel is that, in general, it cannot be integrated over $M$ to produce a formula for the trace of the heat kernel because the coefficient functions are not always integrable over $M$.  In comparison, relatively little corresponding work has been undertaken for the wave kernel on a Riemannian foliation, other than the work of Y. Kordyukov, \cite{Kor}.  In that paper, the researcher derives a trace formula for positive self-adjoint transversally elliptic operators whose principal symbols satisfy certain invariance properties with respect to the leaves of the foliation, using techniques from non-commutative geometry.  These techniques involve representing operators by smooth compactly-supported kernels on the holonomy groupoid of the foliation.  However, this case does not apply to the case of the basic Laplacian, due to the complex nature of the basic projection operator, which cannot generally be represented by such kernels. 
%
%

Recall that for the ordinary Laplacian, the singularities of the trace of the wave kernel contain many spectral invariants; in particular, the spectral invariants associated to the heat kernel can be obtained from the singularity of the trace of the wave kernel at $t=0$. It is natural to examine the possibility of computing additional invariants of the basic spectrum by considering the trace of the basic wave kernel at $t\not=0$. This study of the basic wave kernel has apparently never been undertaken, so the application of wave trace and microlocal techniques to this setting, although straightforward, appears to be new.

This problem is interesting from several points of view.  Riemannian foliations are of interest, both from a geometric point of view as a generalization of a space that is locally a product of Euclidean spaces, and also as a setting for problems in mathematical physics.  (See, for example, the introduction of \cite{GK}.)  In addition, this topic is an extension of the microlocal point of view to the setting of foliations:  here one seeks to make connections between the global geometry of a Riemannian foliation and analysis by associating geometric objects on $Q^*=T^*M/T^*\mathcal{F}\cong N^*\mathcal{F}$ with analytic objects like operators and generalized functions, that are in some sense ``basic".  This particular problem is also of interest from the point of view of spectral theory since the setting allows us to study the spectrum of an operator with a large kernel, and to associate properties of the spectrum with the global geometric structure of the foliation.  Finally, one can regard this problem as being related to the spectral analysis on the space of leaf closures, which is generally quite singular as a space.

The paper is organized as follows: In Section 2, we examine the setting and hypotheses in more detail and define terminology.  We also extend the notion of basic functions to basic generalized functions, and establish the existence of the basic wave kernel. In Section 3, we present the main results about the trace of the basic wave kernel. Section 4 contains the proofs of these results, and Section 5 contains examples of non-simple foliations defined by suspensions, which illustrate the results. 

For general background and notation on Riemannian foliations, see, for example, \cite{Molino}, \cite{T1}, \cite{T2}.

Acknowledgements:  I would like to thank Alejandro Uribe for several helpful conversations.

\end{section}
%
%
\begin{section}{The Setting and Basic Results}

In this section, we carefully examine the setting and investigate its structure.  We then interpret the hypotheses in relation to this structure. After extending the notions of basic functions to generalized functions, we present some elementary results, and define the basic wave kernel on functions.  In particular, we demonstrate that the basic wave kernel exists when the mean curvature form is basic. 

\begin{subsection}{The Stratification of $(M,\mathcal{F})$ and Holonomy}

To study the underlying space of leaf closures, $M/\overline{\mathcal{F}},$ we examine the basic functions. As previously noted in the introduction, the leaves of an arbitrary non-simple foliation are not closed, although the closure of any leaf is a union of leaves, which is an embedded submanifold of $M$.  In fact, each leaf closure is foliated by the leaves that it contains. In general, the leaves have closures of variable dimension, and are defined by a variable dimensional completely integrable distribution $T\overline{\mathcal{F}}$. Furthermore, there exists a natural stratification of $M$ (Section 5.4, \cite{Molino}) as follows.  Let $d(x)$ be the function that assigns to a point $x\in M$ the dimension of the leaf closure containing $x.$  This function takes its values in the positive integers $\{p+k\}$ where $k$ ranges over $0\le k_{1}\le k\le k_{N}\le q,$ with $k_1$ and $k_N$ denoting the minimal and maximal values for $k$, respectively. From \cite{Molino}, Chapter 5, it is known that the function $d(x)$ is lower semi-continuous on $M.$  Let $\Sigma_{p+k}$ denote the inverse image $\{d^{-1} (p+k)\}$. Each of such $\Sigma_{p+k}$ is the (possibly disconnected) union of leaf closures of dimension $p+k$, and is, in fact, an embedded manifold, referred to as the stratum of dimension $p+k$. Furthermore, each stratum $\Sigma_{p+k}$ is foliated by the $p+k$ dimensional leaf closures, by Lemma 5.3 of \cite{Molino}. The lower semi-continuity of $d(x)$ implies the stratum for which the dimension of the leaf closures is maximal is an open dense set in $M$, known as the regular stratum, denoted by $\Sigma_{p+k_N}=\Sigma_{max}.$ (The assumption that there are no dense leaves implies that $k_N<q$.) The strata for which the leaf closures are not of maximal dimension are often referred to all together as the singular strata.   The lower semi-continuity of $d(x)$ implies that for each stratum $\overline{\Sigma_{p+k}}\subset \cup_{\ell\le k}\Sigma_{p+\ell}.$ This partition of $M$ into leaf closures of variable dimension is an example of a singular Riemannian foliation.  

Recall from \cite{W} the holonomy groupoid $\mathcal{G}(\mathcal{F}).$  It has the structure of manifold of dimension $n+p$.  Its elements $\boldsymbol{\alpha}$ are ordered triples $\boldsymbol{\alpha}=\bigl[x,y,[\alpha] \bigr]$ where $x$ and $y$ are points belonging to the same leaf $L$ of $(M,\mathcal{F})$ and $[\alpha]$ is an equivalence class of piecewise smooth curves lying entirely in $L$ with $x=\alpha(0)$ and $y=\alpha(1).$ Its elements define local deffeomorphisms $h_\alpha$ of local transversals in the usual fashion, and via the infinitesimal holonomy map $dh_\alpha,$ define a holonomy action on certain transverse covectors as follows.  Let $V$ be a distribution in $TM=T\mathcal{F}\oplus N\mathcal{F}$. One defines the space of covectors that are transverse to $V$ as follows:  let the subspace $V^0\subset T^*M\setminus \{0\}$ be given by $V^0=\{\xi_x\,|\,\forall X\in V_{x}, i_X(\xi)=0\}.$ We will be interested in the space of covectors that are transverse to the foliation: $(T\mathcal{F})^0.$ The natural action of $\mathcal{G}(\mathcal{F})$ on $(T\mathcal{F})^0$ to is defined for $\xi_{x}\in (T_{x}\mathcal{F})^0$ by
\begin{equation}\label{e:holaction}
\forall X_{y}\in N_y\mathcal{F}\quad (\boldsymbol{\gamma}\cdot\xi)_{y}(X_{y})=\xi_{x}(dh_\gamma^{-1}(X_{y})),
\end{equation}
where $dh_\gamma:N_{x}\mathcal{F}\rightarrow N_{y}\mathcal{F}$ is the differential of the holonomy map of the holonomy element $\boldsymbol{\gamma}.$ A function $F(\xi_x)$ on $T^*M$ will said to be holonomy invariant if $F(\boldsymbol{\gamma}\cdot \xi_y)=F(\xi_x)$ for all $\alpha\in \mathcal{G}(\mathcal{F}).$  

Later, we will need a similar notion of holonomy for the leaf closures for the leaf closures contained in the regular stratum. Observe that one can similarly define the holonomy groupoid $\mathcal{G}_k,$ for each stratum $\Sigma_{p+k},$ where now we simply substitute $\bar{L}\subset \Sigma_{p+k}$ for $L$ and $T\overline{\mathcal{F}}$ for $T\mathcal{F}$ in the discussion above where each $\mathcal{G}_k=\cup_{\overline{L}\subset \Sigma_{p+k}}\mathcal{G}(\mathcal{F})$ is the holonomy groupoid associated to the foliation $(\Sigma_{p+k},\overline{\mathcal{F}}).$ Note that $\mathcal{G}_k$ acts on (suitable) transverse covectors in $N^*\overline{\mathcal{F}}_k:=(T\overline{\mathcal{F}}_k)^0\subset T^*\Sigma_{p+k}.$  Thus, for $k=k_N,$ $\mathcal{G}_{k_N}$ acts on (suitable) transverse covectors in $(T\overline{\mathcal{F}}_{max})^0\subset T^*\Sigma_{max}=T^*M|_{\Sigma_{max}}$.  For $k<k_N$, there is no such action on any transverse covectors in $T^*M|_{\Sigma_{p+k}}$, other than that given by $\mathcal{G}(\mathcal{F})$ holonomy.

\end{subsection}
\begin{subsection}{The Symplectic Setting}

For the purposes of performing microlocal analysis, we will be interested in the symplectic interpretation with respect to $T^*M$ of the various geometric assumptions and structures associated to a foliated manifold with a bundle-like metric.  

For the moment, we will place no conditions on the metric.  The splitting $TM=T\mathcal{F}\oplus T\mathcal{F} ^\perp,$ given by the foliation induces a splitting of $T^*M=(T\mathcal{F})^0\oplus (N\mathcal{F})^0$ where $(T\mathcal{F})^0=\{\xi_x\in T^*M\,|\, \forall X_x\in T_x\mathcal{F}, \,\xi_x(X_x)=0\},$ and is naturally identified $(TM/T\mathcal{F})^*$ and $(NF)^0$ is defined similarly.  Let $\xi_x=(\xi'_x,\xi''_x)$ denote a decomposition of $\xi_x$ with respect to the splitting. Now consider $H(\xi_x)=|\xi_x|_x^2.$  The splitting above of $T^*M$ implies that the function $H$ splits, by the Pythagorean Theorem into $H(\xi_x)=|\xi|^2_x=H_{\mathcal{F}}(\xi_x')+H_{\mathcal{F}^\perp}(\xi_x'').$  The condition that the metric on $M$ should be bundle-like implies that there exists a $\mathcal{G}(\mathcal{F})$ holonomy invariant function $G:(T\mathcal{F})^0\rightarrow \mathbb{R}^+$
such that with respect to the functions $H(\xi_x),$ $H_{\mathcal{F}}(\xi_x)$, and $H_{\mathcal{F}^\perp}(\xi_x)$ the following hold:  (1) $(TF)^0=\{\xi_x\in T^*M\,|\, H_\mathcal{F}(\xi_x)=0\}$;
(2) $H_{\mathcal{F}^\perp}(\xi_x)|_{(TF)^0}=G(\xi_x);$ (3) the hamiltonian vector field $X_\mathcal{F}$ of $H_{\mathcal{F}}$ is zero at all points in $(T\mathcal{F})^0;$ and (4) the hamiltonian vector field $X_{\mathcal{F}^\perp}$ of $H_{\mathcal{F}^\perp},$ satisfies $X_{\mathcal{F}^\perp}=X_H$ and is tangent to $(T\mathcal{F})^0.$  The function $G$ is essentially just the symbol of $\Delta_B.$ In the usual interpretation in terms of the transverse metric $g_T$ on $TM,$ $G(\xi_x)=g^T(\xi_x,\xi_x)$ where $g^T$ denotes the transverse metric induced by the bundle-like metric on $M$ on $N^*\mathcal{F}.$ Note that the metric on $M$ implies that $N^*\mathcal{F},$ the dual of $N\mathcal{F}$ can be identified with $(T\mathcal{F})^0.$  

Notice that $N^*\mathcal{F}$ is a coisotropic submanifold of $T^*M,$ with respect to the usual symplectic form $\omega.$  As such $N^*\mathcal{F}$ is itself foliated by the directions in which the pull-back by the inclusion map $\iota:N^*\mathcal{F}\rightarrow T^*M$ vanishes, that is--the null foliation, which we will denote by $(N^*\mathcal{F},\widetilde{\mathcal{F}})$.  The distribution defining this foliation is precisely the distribution $T\widetilde{\mathcal{F}},$ defined by the canonical lifts of the vector fields $X\in T\mathcal{F}$ to $T^*M,$ which belong to the kernel of $\iota^*\omega.$ The leaves of this foliation through $\xi_x\not=0$ are  
\begin{equation}
\mathcal{L}_{\xi_x}=\{\eta_y\in N^*\mathcal{F}\,|\, \eta_y=(dh^{-1}_\alpha)^*(\xi_x) \exists\,\boldsymbol{\alpha}\in \mathcal{G}(\mathcal{F}), \alpha(0)=x, \alpha(1)=y\}
\end{equation}
Furthermore, the function $G$ is constant along the kernel of $\iota^*\omega$ by $\mathcal{G}(\mathcal{F})$-holonomy invariance:
\begin{equation}
dG_{\xi_x}(X)=0 \text{ for all }X\in Ker(\iota^*\omega), \xi_x\in N^*\mathcal{F}\setminus\{0\}.
\end{equation}
Henceforward, we will delete the zero section from all symplectic manifolds and submanifolds under consideration.
Adopting the notation of \cite{Kor}, let $N\widetilde{\mathcal{F}}$ denote the transverse distribution to the foliation $(N^*\mathcal{F},\widetilde{\mathcal{F}})$--with $\\T(N^*\mathcal{F})=T\widetilde{\mathcal{F}}\oplus N\widetilde{\mathcal{F}}.$  Let $H:T(N^*\mathcal{F})\rightarrow N\widetilde{\mathcal{F}}$ denote the horizontal projection.  The action of holonomy on points in $N^*\mathcal{F}$ induces a lifted holonomy action on $N\widetilde{\mathcal{F}}$ as follows:  for any $\alpha\in \mathcal{G}(\mathcal{F})$ such that $\alpha(0)=x$ and $\alpha(1)=y$ and $\xi_x$ with $\eta_y=(dh_\alpha^{-1})^*(\xi_x)$:
\begin{equation}\label{e:liftedhol}
d\widetilde{h}_{(\alpha,\xi_x)}: N_{\xi_x}\widetilde{\mathcal{F}}\rightarrow N_{\eta_y}\widetilde{\mathcal{F}}.
\end{equation}

%

%
Let $\Phi^t(x,\xi)$ denote the hamiltonian curve associated to $H(x,\xi)^{1/2}=|\xi_x|_x^2.$ Conditions (2) and (4) above imply that the $\Phi^t(x,\xi)$ restricts to $N^*\mathcal{F}$, where $H(x,\xi)^{1/2}=H_{\mathcal{F}^\perp}^{1/2},$ where $\Xi |N^*\mathcal{F}=X_{\mathcal{F}^\perp}.$  An important property of the hamiltonian flow is the following:
\begin{lem}\label{l:phlow}
The transverse flow preserves the leaves of the null-foliation:
$\Phi^t(\mathcal{L}_{\xi_x})=\mathcal{L}_{\Phi^t(\xi_x)}.$  Furthermore, with respect to the splitting of 
$T(N^*\mathcal{F})=T\widetilde{\mathcal{F}}\oplus N\widetilde{\mathcal{F}}$ the differential of the flow splits
\begin{eqnarray}
d\Phi^t & : & T_{\xi_x}\widetilde{\mathcal{F}}\rightarrow T_{\Phi^t(\xi_x)}\widetilde{\mathcal{F}}\\
d\Phi^t & : & N_{\xi_x}\widetilde{\mathcal{F}}\rightarrow N_{\Phi^t(\xi_x)}\widetilde{\mathcal{F}},
\end{eqnarray}
and this last map is preserves $\omega$.
\end{lem}
\begin{proof}
This is a consequence of the fact that $\Phi^t$ is a symplectic diffeomorphism and that the foliation of $N^*\mathcal{F}$ by $T\widetilde{\mathcal{F}}$ is the null-foliation:  $\Phi^t$ preserves the kernel of $\iota^*\omega.$
\end{proof}
The discussion of the previous paragraphs applies equally well to the connected components of the regular stratum $\Sigma_{max}.$  Recall that the inclusion map $\iota:\Sigma_{max}\rightarrow M$ is an embedding, and also that the transverse space of covectors $(T\overline{\mathcal{F}}|_{\Sigma_{max}})^0$ is a subspace of $T^*\Sigma_{max}=T^*M|_{\Sigma_{max}}$ on which $\mathcal{G}_{max}$ acts. Furthermore, the bundle-like metric on $M$ induces a metric which is bundle-like on $\Sigma_{max},$ \cite{Molino}, Chapter 5.4. The discussion of the previous paragraphs implies that $N^*\overline{\mathcal{F}}_{max}$ is a submanifold of the symplectic manifold $T^*M|_{\Sigma_{max}}$ with respect to the usual symplectic form $\omega.$  As above, $N^*\overline{\mathcal{F}}_{max}$ is a coisotropic submanifold of $T^*M|_{\Sigma_{max}}$, and, as such, admits a foliation defined by $T\widetilde{\mathcal{F}}_{max}$, the canonical lift of $T\overline{\mathcal{F}}_{max}$ to $T^*M|_{\Sigma_{max}}.$ Let $N\widetilde{\mathcal{F}}_{max}$ be the corresponding transverse space, and let $H_{max}:T(N^*\mathcal{F})\rightarrow N\widetilde{\mathcal{F}}_{max}$ denote the corresponding horizontal projection. The (non-zero) leaves of this foliation through $\xi_x\in N\widetilde{\mathcal{F}}_{max}$ are given by
\begin{equation}
\mathcal{L}^{max}_{\xi_x}=\{\eta_y\in N^*\overline{\mathcal{F}}_{max}\,|\, \eta_y=(dh^{-1}_{\overline{\alpha}})^*(\xi_x) \exists \boldsymbol{\overline{\alpha}}\in \mathcal{G}_{max}, \overline{\alpha}(0)=x, \overline{\alpha}(1)=y\}.
\end{equation}

The leaves $\mathcal{L}^{max}_{\xi_x}$ are related to the null leaves of $N^*\mathcal{F}$ as follows:
\begin{lem}\label{l:holred}
For $\xi_x\in N^*\overline{L}\subset N\widetilde{\mathcal{F}}_{max},$ each $\mathcal{L}^{max}_{\xi_x}$ is saturated by the leaves of the null-foliation.
\end{lem}
\begin{proof}
Let $\eta_y\in\mathcal{L}_{\xi_x}.$ Hence, there exists a holonomy element $\boldsymbol{\alpha}\in\mathcal{G}(\mathcal{F})$ such that $\alpha(0)=x$ and $\alpha(1)=y$ and $\eta_y=(dh_\alpha^{-1})^*\xi_x.$  Let $\overline{\mathcal{T}}(x)$ a local transversal at $x$ for $(\Sigma_{max},\overline{\mathcal{F}})$ and $\mathcal{T}(x)$ a local transversal for $(M,\mathcal{F})$ containing $\overline{\mathcal{T}}(x),$ and similarly for $\mathcal{T}(y).$ It is sufficient to show that there exists some $\boldsymbol{\overline{\gamma}}\in \mathcal{G}_{max}$ such that $\eta_y=(dh_{\overline{\gamma}}^{-1})^*\xi_x.$ But this is a consequence of the fact that for $\boldsymbol{\alpha}\in\mathcal{G}(\mathcal{F}),$ $\alpha$ also represents a holonomy element $\overline{\alpha}\in\mathcal{G}_{max}$ with $\overline{\alpha}(0)=x$ and $\overline{\alpha}(1)=y.$ We then have the following commutative diagram:
\begin{equation}
\xymatrix{\overline{\mathcal{T}}(x)\ar[r]^{h_{\overline{\alpha}}}\ar[d]_{\iota}&\overline{\mathcal{T}}(y)\ar[d]^{\iota}\\ \mathcal{T}(x)\ar[r]^{h_{\alpha}}&\mathcal{T}(y)
}
\end{equation}
This yields the following commutative diagram:
\begin{equation}
\xymatrix{N_{x}\overline{\mathcal{F}}_{max} \ar[r]^{dh_{\overline{\alpha}}}\ar[d]_{d\iota}&N_y\overline{\mathcal{F}}_{max}\ar[d]^{d\iota}\\ N_x\mathcal{F}\ar[r]^{dh_{\alpha}}&N_y\mathcal{F}
}
\end{equation}
This induces the following commutative diagram
\begin{equation}
\xymatrix{N_x^*\mathcal{F} \ar[r]^{(dh_{\alpha}^{-1})^*}\ar[d]_{d\iota^*}&N^*_y\mathcal{F}\ar[d]^{d\iota^*}\\ N^*_{x}\overline{\mathcal{F}}_{max}\ar[r]^{(dh_{\overline{\alpha}}^{-1})^*}&N^*_y\overline{\mathcal{F}}_{max}
}
\end{equation}
Hence, $(d\iota)^*(dh_{{\alpha}}^{-1})^*=(dh_{\overline{\alpha}}^{-1})^*(d\iota)^*$ for all $\alpha\in\mathcal{G}(\mathcal{F})$. But $d\iota$ is just the identity map on $N_{x}\overline{\mathcal{F}}_{max}$, and $d\iota^*$ is just the orthogonal projection onto $N_{x}^*\overline{\mathcal{F}}_{max}.$ We conclude that
$(dh_\alpha^{-1})$ restricts to the image of $d\iota,$ which is just $N_{x}\overline{\mathcal{F}}_{max}$ and thus
\begin{equation}
\eta_y=(dh_\alpha^{-1})^*(\xi_x)=(dh_{\overline{\alpha}}^{-1})^*d\iota^*(\xi_x),
\end{equation}
which proves the result.
\end{proof}

In a similar fashion to \eqref{e:liftedhol}, the action of holonomy on points in $N^*\overline{\mathcal{F}}_{max}$ induces a lifted holonomy action on $N\widetilde{\mathcal{F}}_{max}$ as follows:  for any $\boldsymbol{\overline{\alpha}}\in \mathcal{G}_{max}$ such that $\overline{\alpha}(0)=x$ and $\overline{\alpha}(1)_y$ and $d\iota^*(\xi_x)$ with $d\iota^*(\eta_y)=(dh_{\overline{\alpha}}^{-1})^*(d\iota)^*(\xi_x)$:
\begin{equation}\label{e:liftedholk}
d\widetilde{h}^{max}_{(\overline{\alpha},\xi_x)}: N_{\xi_x}\widetilde{\mathcal{F}}_{max}\rightarrow N_{\eta_y}\widetilde{\mathcal{F}}_{max}.
\end{equation}
Finally, as a corollary to the reasoning of the previous paragraphs and of Lemma \ref{l:phlow}, we have the analogous result:
\begin{lem}\label{l:phlowbar}
The transverse flow restricts to $N^*\widetilde{\mathcal{F}}_{max}$ over the connected components of $\Sigma_{max}$, and preserves the leaves of the foliation by $\mathcal{L}^{max}_{\xi_x}:$ 
$\Phi^t(\mathcal{L}^{max}_{\xi_x})=\mathcal{L}^{max}_{\Phi^t(\xi_x)}.$
With respect to the splitting of 
$T(N^*\overline{\mathcal{F}}_{max})=T\widetilde{\mathcal{F}}_{max}\oplus N\widetilde{\mathcal{F}}_{max}$ the differential of the flow splits
\begin{eqnarray}
d\Phi^t & : & T_{\xi_x}\widetilde{\mathcal{F}}_{max}\rightarrow T_{\Phi^t(\xi_x)}\widetilde{\mathcal{F}}_{max}\\
d\Phi^t & : & N_{\xi_x}\widetilde{\mathcal{F}}_{max}\rightarrow N_{\Phi^t(\xi_x)}\widetilde{\mathcal{F}}_{max}.
\end{eqnarray}

\end{lem}
\end{subsection}
%
%
\begin{subsection}{Basic Distributions and the Basic Wave Kernel}

In what follows we assume the following for $(M,\mathcal{F}):$ (1) $M$ is equipped with a metric that is bundle-like with respect to the foliation $\mathcal{F} $;  (2) that the foliation is transversally orientable (see below); (3) that the mean curvature form for the foliation is a basic one-form; and (4) the maximal leaf closure satisfies $k_N<q$, and thus, there are no dense leaf closures.

From \cite{PaRi}, there is a natural projection $P$ from $C^\infty(M)\rightarrow C^\infty_B(M,\mathcal{F}).$  In fact, this projection extends to a projection from $L^2(M)\rightarrow L^2_B(M,\mathcal{F}),$ the space of basic functions in $L^2(M).$  The projection $P$  is self-adjoint.  Indeed, we can describe the operator $P$ in terms of a series of push-forwards and pull-backs by submersions as follows:  Let $\widehat{\pi}:\widehat{M}\longrightarrow M$ be the oriented transverse frame bundle. (Note: the assumption that the foliation is transversally orientable is made purely for the sake of simplicity.  If the foliation is not transversally oriented,  $\widehat{M}$ has two connected components, and we replace $\widehat{M}$ by one of these components.) The foliation of $M$ lifts to a $p$ dimensional foliation of $\widehat{M},$  denoted by $(\widehat{M}, \widehat{\mathcal{F}}). $ Let $K$ be a typical leaf in the lifted foliation, $(\widehat{M},\widehat{\mathcal{F}}),$ and let $\bar{K}$ denote the closure of this leaf.  In fact, the closures of the leaves of the lifted foliation, $\bar{K},$ are the fibres of a fibre bundle over a compact manifold, $W,$ called the basic manifold.  Let $\rho:\widehat{M}\rightarrow W$ denote this bundle projection.  We then have the following double fibration of $\widehat{M}$:
\begin{equation}\label{e:doublefib}
\xymatrix{
 &\widehat{M} \ar[dl] _{\hat{\pi}}\ar[dr]^\rho \\
M  &  & W}
\end{equation}
where each of the fibrations is locally trivial, \cite{Molino}.
Let $f\in C^\infty(\widehat{M})$.  Define an operator $A:C^\infty(\widehat{M})\rightarrow C^\infty_B(\widehat{M})$ to be the operator obtained by averaging $f$ over the closures of the leaves of $(\widehat{M},\widehat{\mathcal{F}})$. It has been shown in \cite{PaRi} that $A$ as defined above is formally self-adjoint with respect to the $L^2$ inner product of functions on $\widehat{M}$.  The basic projector can be expressed in terms of these operations as $P=\widehat{\pi}_*A\widehat{\pi}^*.$ 

We can extend this projector to distributions as follows:
\begin{lem}\label{l:basicd}
Let $\langle \cdot,\,\cdot\rangle_{X}$ denote the pairing of a generalized function with a function on a manifold $X$. Let $u\in\mathcal{D'}(M)$ and $\nu\in\mathcal{D}'(\widehat{M}).$  We extend $A$ to generalized functions, by $\langle A\nu,\varphi\rangle_{\widehat{M}}=\langle \nu,A\varphi\rangle_{\widehat{M}},$ where $\varphi\in C^\infty(\widehat{M}).$  We similarly extend $P$ to generalized functions similarly: $\langle Pu,\varphi\rangle_{M}=\langle u,P\varphi\rangle_{M},$ for $\varphi\in C^\infty(M)$.
 \end{lem}
 
\begin{proof}
Since $\widehat{\pi}$ is a submersion, $\widehat{\pi}^*u$ and $\widehat{\pi}_*\nu$ are a well-defined generalized functions on $\widehat{M}$ and $M,$ respectively.  Now consider the operator $A$. From Lemma 1.4 of \cite{PaRi}, for all $f,\, g\in L^2(\widehat{M}),$ $\langle Ag,f\rangle_{L^2}=\langle g,Af\rangle_{L^2},$ where $\langle\cdot,\cdot\rangle_{L^2}$ is the $L^2(\widehat{M})$ inner product.  Since this inner product coincides with the pairing of a generalized function with functions on $L^2,$ the corresponding statement $\langle A\nu,f\rangle_{\widehat{M}}=\langle \nu,Af\rangle_{\widehat{M}}$ holds for $\nu\in \mathcal{D}'(\widehat{M}).$
Thus, $A$ is a well-defined operation on any generalized function $\nu\in \mathcal{D}'(\widehat{M}).$ By the usual functorial relations for generalized functions, $P$ is a composition of well-defined operations on generalized functions.
\end{proof}

We next define a notion of what it means for a generalized function to be basic.

{\bf Definition 1:}   A generalized function $u\in\mathcal{D}'(M)$ is basic if $Pu=u,$ or, equivalently, if $\langle u,\varphi \rangle_M=\langle u,P\varphi\rangle_M,$ for all $\varphi\in C^\infty(M).$

\begin{prop}
The following are equivalent:
\begin{itemize}
\item[(1)] The generalized function $u\in\mathcal{D}'(M)$ is basic.
\item[(2)] $X(u)=0$ for every vector field $X\in T\overline{\mathcal{F}}$ defined on $U$ which contains $supp(u)$.
\end{itemize}
\end{prop}

\begin{rem}
If $u$ is a basic generalized function such that $u\notin C^\infty(M),$ then the wave front set of $u$, $WF(u)\subset N^*\mathcal{F}\setminus\{0\}.$ In fact, we will see that the wave front set of $u$ will actually be contained in $\cup_{L\subset M}N^*(\overline{L}).$ Let $N^*\overline{\mathcal{F}}$ denote this set.
\end{rem}

We define the basic wave kernel on functions in an analogous manner to the basic heat kernel:

{\bf Definition 2:} Let $(x,y)$ be coordinates on $M\times M.$  Let $D_t=\frac{1}{i}\frac{\partial}{\partial t}.$ Define the basic wave kernel (acting on functions), $U_B(t,x,y),$ as the solution to the system:
\begin{eqnarray}
(D_t+\sqrt{\Delta_B}_{x})U_B(t,x,y)&=&0 \nonumber\\
U_B(0,x,y)&=&\delta(x-y)\text{ on basic functions.}
\end{eqnarray}
(Here $\sqrt{\Delta_B}$ is can be defined via $\sqrt{\Delta}$ and $P$.)

\begin{rem}
Note that $U_B(t,x,y)$ is generalized function on $\mathbb{R}\times M\times M$ that is basic on each $M$ factor. This is analogous to the basic heat kernel, which is a basic function on each factor, see, for example, \cite{KR2}.
\end{rem}

\begin{thm}
The basic wave kernel $U_B(t,x,y)$ exists.  It is unique solution to
\begin{equation}\label{e:basicwave}
U_B(t,x,y)=P_{x}P_{y}U(t,x,y)=\sum_{j=1}^\infty e^{-it\sqrt{\lambda_j^B}}e_j(x)e_j(y),
\end{equation}
where $U(t,x,y)$ is the wave kernel for the ordinary Laplacian on $M$ and \\
$0\le \lambda_1^B\le \lambda_2^B\le \dots$ are the eigenvalues of $\Delta_B.$ In the above, $P_{x}$ denotes the basic projector acting on the first $M$ factor of $M\times M\times \mathbb{R},$ the space on which the wave kernel is defined.  $P_{y}$ denotes analogously the basic projector on the second factor.
\end{thm}
\begin{proof}
The wave kernel for the ordinary Laplacian is the unique solution to the system
\begin{eqnarray}
(D_t+\sqrt{\Delta}_{x})U(t,x,y)&=&0 \label{e:wave}\\
U(0,x,y)&=&\delta(x-y).
\end{eqnarray}

If we apply the operator $P_{x}P_{y}$ to both sides of \eqref{e:wave}, we have:

\begin{equation}\label{e:basicexist}
(D_tP_{x}P_{y}+P_{x}\sqrt{\Delta}_{x}P_{y})U(t,x,y)=0.
\end{equation}
If we apply Theorem 2.7 of \cite{PaRi} specialized to a foliation with basic mean curvature,
we have that $\Delta_BP=P\Delta.$  It follows that $\sqrt{\Delta_B}P=P\sqrt{\Delta}.$ Hence, equation \eqref{e:basicexist} becomes
\begin{equation}
(D_tP_{x}P_{y}+\sqrt{\Delta_B}_{x}P_{x}P_{y})U(t,x,y)=0,
\end{equation}
thus proving \eqref{e:basicwave}.  The  initial condition is immediate.
\end{proof}
We wish to compute the basic wave trace--that is, if $\Pi:\mathbb{R}\times M\rightarrow \mathbb{R}$ and $\Delta:M\rightarrow M\times M$ is the diagonal map, then we wish to compute:
\begin{equation}
\Pi_*\Delta^*U_B(t,x,y)=\Pi_*\Delta^*P_{x}P_{y}U(t,x,y).
\end{equation}
Note that  
\begin{equation}
Trace \,U_B(t)=\Pi_*\Delta^*U_B(t,x,y)=\sum_{j=1}^\infty e^{-it\sqrt{\lambda^B_j}},
\end{equation}
which is just the Fourier transform of the spectral distribution of $\sqrt{\Delta_B}:$  $$\sigma(\mu)=\sum_{j=1}^\infty\delta(\mu-\sqrt{\lambda^B_j}).$$

\end{subsection}
\begin{subsection}{Relatively Closed Curves with Respect to $(M,\overline{\mathcal{F}})$}
Of particular interest in this analysis are certain arcs of the hamiltonian curves of the transverse metric:

{\bf Definition 3:}  An arc of a curve $\Phi^T$ in $N^*\mathcal{F}$ will said to be relatively closed with respect to the (singular) foliation $\overline{\mathcal{F}}$ with relative period $T$ if its endpoints $\xi_x$ and $\eta_y=\Phi^T(\xi_x)$ belong to $N^*\overline{L}$ for $L\subset \Sigma_{p+k}$ and either (1)  $k<k_N$  and $\Phi^T(\xi_x)\in \mathcal{L}_{\xi_x}$; or  (2) $k=k_N$ and $\Phi^T(\xi_x)\in \mathcal{L}^{max}_{\xi_x}.$ The projection $\gamma(t,x)=\pi(\Phi^t(x,\xi))$ of a relatively closed hamiltonian curve in $(N^*\mathcal{F},\widetilde{\mathcal{F}})$ by $\pi: N^*\mathcal{F}\rightarrow M$ will be said to be {\it relatively closed} with respect to the singular foliation $(M,\overline{\mathcal{F}}).$ 

Note that in local distinguished coordinates, it is easily seen that $\gamma'(0,x)=(d\pi)_{\xi_x}(\Xi)\perp T_x\overline{\mathcal{F}}.$ From Chapter 6 of \cite{Molino}, it is known that if $\gamma(t,x)$ is a geodesic passing through $x$, that is perpendicular to the leaf closures, then it remains perpendicular to all the leaf closures that it meets. Thus, the projections of such relatively closed hamiltonian curves are geodesic arcs that are orthogonal to the leaf closures through which the geodesic passes.

Now consider the set of endpoints of relatively closed of the hamiltonian flow $\Phi^T$ restricted to $N^*\mathcal{F}$:  $Z_T=\bigcup_kZ_T^k,$ where for each $k<k_N$ each $Z^k_T$ is given by
\begin{equation}\label{e:zk}
Z^k_T=\{\xi_x\in N_x^*\overline{L},\, \overline{L}\subset \Sigma_{p+k} \,|\, \Phi^T(\xi_x)\in\mathcal{L}_{\xi_x}\}
\end{equation}
and for $k=k_N$, $Z^{k_N}_T=Z_T^{max}$ is given by
\begin{equation}\label{e:zmax}
Z^{max}_T= \{\xi_x\in N_x^*\overline{L},\,L\subset \Sigma_{max}\,| \,\Phi^T(\xi_x)\in \mathcal{L}^{max}_{\xi_x}\}.
\end{equation}
Note that $Z^{max}_T\cap S(N^*\mathcal{F})$ and $Z^0_T\cap S(N^*\mathcal{F})$ are closed.

The set $Z_T$ is also a saturated by the null-leaves.  First note the following 
\begin{lem}\label{l:holonomy}
Suppose $\iota:S\rightarrow M$ is a saturated embedded submanifold.  Then $N^*S$ is saturated by leaves of the null-foliation. 
\end{lem}
\begin{proof}
Suppose $\bar{x}$, $\bar{y}\in S$ such that $\iota(\bar{x})=x$ and $\iota(\bar{y})=y,$ and let $\mathcal{T}_S(\bar{x})$ and $\mathcal{T}_S(\bar{y})$ be local transversals in $S$ at $\bar{x}$ and $\bar{y}$, respectively.  Let $\mathcal{T}(x)$ and $\mathcal{T}(y)$ be local transversals in $M$ containing $\iota(\mathcal{T}_S(\bar{x}))$ and $\iota(\mathcal{T}_S(\bar{y}))$. Given any holonomy element $\boldsymbol{\alpha}\in\mathcal{G}(\mathcal{F}),$ with $\alpha(0)=x$ and $\alpha(1)=y$, we have a local diffeomorphism $h_\alpha:\mathcal{T}(x)\rightarrow \mathcal{T}(y)$. There is a corresponding holonomy element in $\mathcal{G}(S,\mathcal{F})$, which we also denote by $\boldsymbol{\alpha}$, and the corresponding diffeomorphism will be denoted by $h^S_{\alpha}:\mathcal{T}_S(\bar{x})\rightarrow \mathcal{T}_S(\bar{y}).$  We then have the following commutative diagram:
\begin{equation}
\xymatrix{\mathcal{T}_S(\bar{x})\ar[r]^{h^S_\alpha}\ar[d]_{\iota}&\mathcal{T}_S(\bar{y})\ar[d]^{\iota}\\ \mathcal{T}(x)\ar[r]^{h_\alpha}&\mathcal{T}(y)
}
\end{equation}
This yields the following  commutative diagram:
\begin{equation}
\xymatrix{N_{\bar{x}}\mathcal{F}_S\ar[r]^{dh^S_\alpha}\ar[d]_{d\iota}&N_{\bar{y}}\mathcal{F}_S\ar[d]^{d\iota}\\ N_x\mathcal{F}\ar[r]^{dh_\alpha}&N_y\mathcal{F}
}
\end{equation}
where $N_{\bar{x}}\mathcal{F}_S$ is the transverse space in $T_{\bar{x}}S$.  Let $N(T_x\mathcal{F},T_xS)\subset T_xM$ denote the complement of $T_x\mathcal{F}$ in $T_xS$ regarded as a subspace of $T_xM$. 
Then $d\iota(N_{x}\mathcal{F}_S)=N(T_x\mathcal{F},T_xS)$ which yields the commutative diagram below:
\begin{equation}
\xymatrix{N_{\bar{x}}\mathcal{F}_S\ar[r]^{dh^S_\alpha}\ar[d]_{d\iota}&N_{\bar{y}}\mathcal{F}_S\ar[d]^{d\iota}\\ N(T_x\mathcal{F},T_xS)\ar[r]^{dh_\alpha}&N(T_y\mathcal{F},T_yS)
}
\end{equation}
Thus, the infinitesimal holonomy map restricts: \\
$dh_\alpha:N(T_x\mathcal{F},T_xS)\rightarrow N(T_y\mathcal{F},T_yS)$.  

Now suppose $\xi_x\in N^*S$ and let $\eta_y\in \mathcal{L}_{\xi_x}$ with $\eta_y=(dh_\gamma^{-1})^*\xi_x.$  For all $X_x\in N(T_x\mathcal{F},T_xS),$ $\xi_x(X_x)=0$ by definition of $N^*S$. If $Y_y\in N(T_y\mathcal{F},T_yS),$ then there is an $X_x\in N(T_x\mathcal{F},T_xS)$ with $Y_x=dh_\gamma(X_x)$ because the holonomy action is invertible.  It then follows that $\eta_y(Y_y)=\xi_x(X_x)=0,$ and hence $\eta_y\in N^*S$, and the result follows.
\end{proof}

From this, we see the following.
\begin{lem}\label{l:ZTsat}
The set $Z_T$ is saturated by the leaves of the null foliation, and $Z_T^{max}$ is saturated by the leaves  $\mathcal{L}^{max}_{\xi_x}$.
\end{lem}
\begin{proof}
To prove the first part of the lemma, we need only show that each $Z_T^k$ is saturated.  For $k<k_N$, the fact that $Z_T^k$ is foliated follows from Lemma \ref{l:phlow}:  Let $\xi_x\in Z_T^k$ and let $\eta_y\in \mathcal{L}_{\xi_x}$.  By hypothesis, $\Phi^T(\xi_x)\in \mathcal{L}_{\xi_x},$ so 
\begin{equation}
\eta_y\in\Phi^T(\mathcal{L}_{\xi_x}=\mathcal{L}_{\Phi^T(\xi_x)}=\mathcal{L}_{\Phi^T(\eta_y)}.
\end{equation}
For $k=k_N$, suppose $\xi_x\in N^*\overline{L}$ and consider $\eta_y\in \mathcal{L}_{\xi_x}$ with holonomy element $\beta$ such that $\beta(0)=x$ and $\beta(1)=y$ such that $\eta_y=(dh^{-1}_\beta)^*\xi_x.$  By the previous lemma, such a covector $\eta_y$ is also in $N^*\overline{L}$ by the previous lemma, and so $\eta_y=d\iota^*\eta_y.$ As in the proof of Lemma \ref{l:holred}, each $\beta\in \mathcal{G}(\mathcal{F})$ represents a $\mathcal{G}_{max}$ holonomy element $\overline{\beta}$ with $\overline{\beta}(0)=x$ and $\overline{\beta}(1)=y$  Then, using the fact that $d\iota^*(dh_\beta^{-1})^*=(dh_{\overline{\beta}}^{-1})^*d\iota^*$ from the proof of Lemma \ref{l:holred}, it follows that 
\begin{eqnarray*}
\eta_y&=&d\iota^*(\eta_y)=d\iota^*(dh_\beta^{-1})^*\xi_x=(dh_{\overline{\beta}}^{-1})^*d\iota^*(\xi_x)\\
&=&(dh_{\overline{\beta}}^{-1})^*(dh_{\overline{\alpha}}^{-1})^*d\iota^*(\Phi^T(\xi_x))
\end{eqnarray*}
where $\overline{\alpha}(0)=\pi(\Phi^T(\xi_x))$ and $\overline{\alpha}(1)=x$ (since $\Phi^T(\xi_x)\in \mathcal{L}^{max}_{\xi_x}$.)Thus $\eta_y=(dh_{\overline{\gamma}_1}^{-1})^*d\iota^*(\Phi^T(\xi_x)),$
where $\overline{\gamma}_1$ is $\overline{\alpha}\circ\overline{\beta}.$ Now recall that $\Phi^T(\mathcal{L}_{\xi_x})=\mathcal{L}_{\Phi^T(\xi_x)},$ so $d\iota^*\Phi^T(\xi_x)=\Phi^T(\xi_x),$ and, furthermore, there is a $\gamma_2\in \mathcal{G}(\mathcal{F})$ with $\gamma_2(0)=\pi(\Phi^T(\xi_x))$ and $\gamma_2(1)=\pi(\Phi^T(\eta_y))$ such that $\Phi^T(\xi_x)=(dh_{\gamma_2}^{-1})^*\Phi^T(\eta_y).$ Then we have 
\begin{eqnarray*}
\eta_y&=&(dh_{\overline{\gamma}_1}^{-1})^*(dh_{\gamma_2}^{-1})^*d\iota^*\Phi^T(\eta_y)\\
&=&(dh_{\overline{\gamma}_1}^{-1})^*(dh_{\overline{\gamma}_2}^{-1})^*d\iota^*\Phi^T(\eta_y)\\
&=&(dh_{\overline{\gamma}}^{-1})^*\Phi^T(\eta_y).
\end{eqnarray*}
where $\overline{\gamma}=\overline{\gamma}_2\circ \overline{\gamma}_1.$ Thus, $\eta_y\in Z_T^k,$ proving the result.  
The second part of the lemma for $Z_T^{max}$ follows by reasoning analogous to that of the first part of the proof, using Lemma \ref{l:phlowbar}.
\end{proof}

The component $Z_T^{max}$ of the relative fixed point set is said to be clean if the following holds.

{\bf Definition 4: }Let $T$ be the length of a relatively closed arc of the hamiltonian flow $\Phi^t.$  We say that the relative fixed point set is $Z_T^{max}$ is clean if (1)$ Z^{max}_T$ is a smooth submanifold of $N^*\mathcal{F};$ and (2) for every $\xi_x\in Z_T^{max}$ with $\eta_y=(dh_{\overline{\alpha}}^{-1})^*\xi_x=\Phi^T(\xi_x)$ then $d\Phi^T_{\xi_x}(T_{\xi_x}Z^{max}_T)=T_{\eta_y}Z^{max}_T$ for $\overline{\alpha}\in \mathcal{G}(\overline{L})$ with $\overline{\alpha}(0)=x$ and $\overline{\alpha}(1)=y$. Note that the condition that $\eta_y=(dh_{\overline{\alpha}}^{-1})^*\xi_x$ implies that for all $\xi_x\in Z_T^{max}$
\begin{equation}
d\Phi^T_{\xi_x}(N_{\xi_x}\widetilde{\mathcal{F}}_{max})=d\widetilde{h}^{max}_{(\overline{\alpha},\xi_x)}(N_{\xi_x}\widetilde{\mathcal{F}}_{max})=N_{\eta_y}\widetilde{\mathcal{F}}_{max}
\end{equation}
by Lemma \ref{l:phlowbar}. Note that the $T\widetilde{\mathcal{F}}$ components of $TZ_T^{max}$ are also determined by \ref{l:phlowbar}.
\begin{rem}
The above definition of clean-ness is with respect to the holonomy of the leaf closures in $\Sigma_{max}.$  One could define a similar notion of clean-ness for relative fixed points in the singular strata, using the $\mathcal{G}(F)$ holonomy as follows:
\end{rem}

{\bf Definition 5: }Let $T$ be the length of a relatively closed arc of the hamiltonian flow $\Phi^t.$  For $k<k_N$, $Z_T^k$ is clean if (1) $Z^k_T$ is a smooth submanifold of $N^*\mathcal{F};$ and (2) for every $\xi_x\in Z_T^k$ with $\eta_y=(dh_{\alpha}^{-1})^*\xi_x=\Phi^T(\xi_x)$ then $d\Phi^T_{\xi_x}(T_{\xi_x}Z^k_T)=T_{\eta_y}Z^k_T$ for $\overline{\alpha}\in \mathcal{G}(\mathcal{F})$ with $\alpha(0)=x$ and $\alpha(1)=y$. Note that the condition that $\eta_y=(dh_\alpha^{-1})^*\xi_x$ implies that for all $\xi_x\in Z_T^k$
\begin{equation}
d\Phi^T_{\xi_x}(N_{\xi_x}\widetilde{\mathcal{F}})=d\widetilde{h}_{(\alpha,\xi_x)}(N_{\xi_x}\widetilde{\mathcal{F}})=N_{\eta_y}\widetilde{\mathcal{F}}.
\end{equation}

\end{subsection}
\end{section}
%
%
\begin{section}{Main Results}

In this section, we present our main results concerning the wave trace of the basic Laplacian.

\begin{thm}\label{t:sojourntimes}
In the notation previously established, 
\begin{eqnarray}\label{e:invariants}
WF\bigl(\Pi_*\Delta^*(U_B(t,x,y) \bigr)\subset\{(T,\tau)\!\!\!\!\!&|&\!\!\!\!\! \tau<0\,,\,T=\text{length of a relatively closed}\nonumber\\ &\quad&\text{geodesic arc.}\}
\end{eqnarray}

\end{thm}

In the notation of \cite{Z}, the lengths of curves such as these are sometimes referred to as ``sojourn times". Let $\mathcal{ST}(M,\mathcal{F})$ denote the set of lengths of relatively closed hamiltonian curves corresponding to the hamiltonian function $H^{1/2}=\sigma(\sqrt{\Delta})$.  These curves project down to relatively closed geodesic arcs for $(M,\mathcal{F}).$  For the sake of convenience, we shall denote the set of lengths of such curves as just $\mathcal{ST}.$  In the analysis that follows, it will be necessary to make a distinction between the relatively closed geodesic arcs that remain inside the regular stratum $\Sigma_{max}$ on $M$ and those that leave the regular stratum.  Let the sojourn times $T$ that correspond to relatively closed  geodesic arcs that remain inside the regular stratum be called {\it regular sojourn times}, and denote the set of such $T$ by $\mathcal{RST}(M,\mathcal{F}).$ Let the sojourn times in the complement of $\mathcal{RST}$ be known as {\it singular sojourn times}, and denote the set of such $T$ by $\mathcal{SST}(M,\mathcal{F}).$  

Now suppose that $\mathcal{SST}\subset U_1$ and $\mathcal{RST}\subset U_2$ where $U_1$ and $U_2$ are disjoint open sets in $\mathbb{R}.$  Then we can pick $\chi\in C^\infty(\mathbb{R})$ with $\chi(t)=0$ on $U_1$ and $\chi(t)=1$ on $U_2.$ In this case, the intersection of the $\mathcal{C}_{max}$ component of the canonical relation of $P$ and $U(t,x,y)$ is clean, and thus, as observed in \cite{Z}, $WF\bigl(\chi(t)\Pi_*\Delta^*(U_B(t,x,y)\bigr)$ is a conic Lagrangian submanifold of $T^*\mathbb{R},$ and thus must be a union of rays over the discrete set of sojourn times $\mathcal{RST}.$ In the notation and terminology of \cite{Z}, we let $\Gamma^T_{max}=\{T,\tau\,|\,\tau<0\}$ denote the ray over $T\in \mathcal{RST}$.  The relatively closed orthogonal geodesic arcs of a given length $T$ make up conic submanifolds $Z^{max}_T$ whose connected components are finite in number and denoted by $Z^{max}_{T,j}.$ Let $S(Z^{max}_{T,j})$ be the set $\{(T,\tau)\in Z^{max}_{T,j}\,|\,|\tau|=1\},$ and let $e_{T,j}:=dim(S(Z^{max}_{T,j}))$ and let $e_T=max\{e_{T,j}\}.$ Finally, we must assume that the set $Z^{max}_T$ of relative fixed points of the hamiltonian flow $\Phi^T$ on $N^*\mathcal{F}$ are clean for all $T\in \mathcal{RST}$ in the sense Definition 4. 


%
\begin{thm}\label{t:partialtrace}
With the above assumptions,  
\begin{equation}
\chi(t)\Pi_*\Delta^*\bigl(U_B(t,x,y)\bigr)=\sum_{T\in \mathcal{RST}}\nu_{T}(t),
\end{equation}
where $\nu_{T}\in I^{-1/4-e_T/2-r}(\mathbb{R},\Gamma^T,\mathbb{R})$ where $r=-(p+k_N)/2.$  Furthermore, $\nu_{T}$ has an expansion of the form
\begin{equation}\label{e:expansion}
\nu_{T}(t)=e^{\frac{i\pi m_T}{4}}\sum_{j=0}^\infty\sigma_j(T)(t-T+i0)^{-\frac{e_T-1}{2}+\frac{p+k_N}{2}-j}\,mod\,C^\infty(\mathbb{R}),
\end{equation}
where $m_T$ is the Maslov index of $Z^{max}_T,$ and $\sigma_j(T)=\int_{S(Z_{T}^{max})}d\mu_{Z_{T,j}^{max}}$ where for each $j\ge 0$, $d\mu_{Z_{T,j}^{max}}$ is a density on $S(Z_{T}^{max})$.
Note also that the rank of $N^*\overline{\mathcal{F}}_{max}$ depends on $k_N$ in the above formula.
\end{thm}


If the foliation has regular closure, then a complete expansion of the trace can be obtained.   Let $k$ be such that all the leaf closures are of dimension $p+k$ ($k:=k_N$).  Since the leaf closures are all of the same dimension, $M$ is foliated by the leaf closures, denoted by $(M,\overline{\mathcal{F}}).$  Let $T\overline{\mathcal{F}}$ be the associated distribution, and let $N\overline{\mathcal{F}}$ be the orthogonal distribution.  In this notation, we have the following:

\begin{cor}\label{t:regclosure}
If $(M,\mathcal{F})$ has regular closure, then 
\begin{equation}
\Pi_*\Delta^*\bigl(U_B(t,x,y)\bigr)=\sum_{T\in \mathcal{ST}}\nu_T(t),
\end{equation}
where $\nu_T\in I^{-1/4-e_T/2-r}(\mathbb{R},\Gamma^T,\mathbb{R})$ where $r=-(p+k)/2$ is the degree of $K\in I^r(M\times M,\mathcal{C},\Omega^{1/2}_{M\times M}).$  Furthermore, $\nu_T$ has an expansion of the form
\begin{equation}
\nu_T(t)=e^{\frac{i\pi m_T}{4}}\sum_{j=0}^\infty\sigma_j(T)(t-T+i0)^{-\frac{e_T-1}{2}+\frac{p+k}{2}-j}\,mod\,C^\infty,
\end{equation}
with the leading term given as in the previous theorem.  Note that under these hypotheses $p+k$ is constant and $r=-(p+k)/2$ is the degree of $K$, which can here be represented as a single Lagrangian distribution $K\in I^r(M\times M,\mathcal{C},\Omega^{1/2}_{M\times M}).$
\end{cor}

If $(M,\mathcal{F})$ does not have regular closure, then clean-ness fails, since then the basic projector $P$ not have a nice canonical relation which is necessarily Lagrangian, and it is not clear if an expansion like the one in \eqref{e:expansion} above exists for {\it all} $T\in \mathcal{ST}$.  

The different nature of the results for the regular closure and more general case is not entirely unexpected.  From the structure theorems for Riemannian foliation (Theorem 5.1, Proposition 5.2 of \cite{Molino}), it is known that if a foliation admits regular closure, then the space of leaf closures $M/\overline{\mathcal{F}}$ has the relatively nice structure of an orbifold. (The cone points in the orbifold structure arise from leaf closures with non-trivial holonomy.)  If the general case, by contrast, the best that can be said about the structure of $M/\overline{\mathcal{F}}$ is that it can be identified with the orbit space of the inherited $SO(q)$ action on the basic manifold $W$.

It is possible that a more complete wave trace result may be available using the additional structure of the double fibration in \eqref{e:doublefib}, in the spirit of \cite{KR2}.  From the theory of Riemannian foliations, we know that the foliation induced by the lifted foliation has regular closure.  However, the mean curvature form associated to the lifted foliation is not necessarily basic, and thus, the basic wave kernel for $\widehat{M}$ may not exist.  Nonetheless, the basic projector for the lifted foliation is the averaging operator $A$, \cite{PaRi},  and by a calculation analogous to the one to follow in Proposition \ref{p:canonical}, it can be shown to be an operator whose Schwartz kernel is a Lagrangian distribution.  A possible approach to the problem of analyzing the singular sojourn times is to make use of this additional structure by representing the basic projector $P$ as $\hat{\pi}_*A\hat{\pi}^*$ and understanding the behavior of the singularities at the non-clean intersection of the canonical relations of $\hat{\pi}^*$ and $A$.  A further use of the double fibration structure would be to relate the relatively closed curves on $M$ with geometrically interesting closed curves on $W$. The double fibration structure yields the map $\rho_*\hat{\pi}^*$ between generalized functions on $M$ and generalized functions on $W,$ that we conjecture yields a correspondence between basic generalized functions on $M$ and some class of generalized functions on $W$ when restricted to the basic generalized functions on $M$. This map appears to be a kind of generalized version of a Radon transform, although in this case the transformation may not be invertible. In this way, we suspect that a more general and satisfying description of the basic wave trace in terms of the geometry of this structure may yet be forthcoming.

Finally, we have the following corollary, in the case that $(M,\mathcal{F})$ has minimal stratum, $\Sigma_p$, of compact leaves of dimension $p$.  Note that by the lower continuity of the leaf closures, this stratum is compact.

Let $\mathcal{MST}$ denote the set of relative periods of hamiltonian curves whose endpoints lie in the minimal stratum $\Sigma_p$. If this set of periods lies in an open set $U_1\subset \mathbb{R}$ and the complementary set of periods $\mathcal{MST}^c\subset U_2$ where $U_1$ and $U_2$ are disjoint open sets, then pick, as before $\chi\in C^\infty(\mathbb{R})$ with $\chi(t)=0$ on $U_1$ and $\chi(t)=1$ on $U_2.$ In this case, the intersection of the component of the canonical relation of $P$ that corresponds to the minimal stratum and $U(t,x,y)$ is clean, and thus, as above, we may cut off the wave trace.  In this case $WF\bigl(\chi(t)\Pi_*\Delta^*(U_B(t,x,y)\bigr)$ will again be a union of rays over the discrete set of sojourn times $\mathcal{MST}.$ If we let $\Gamma^T_0=\{T,\tau\,|\,\tau<0\}$ denote the ray over $T\in \mathcal{MST}$.  The relatively closed orthogonal geodesic arcs of a given length $T$ make up conic submanifolds $Z^0_T$ whose connected components are finite in number and denoted by $Z^0_{T,j}.$ Let $S(Z^0_{T,j})$ be the set $\{(T,\tau)\in Z^0_{T,j}\,|\,|\tau|=1\},$ and let $e_{T,j}:=dim(S(Z^0_{T,j}))$ and let $e_T=max\{e_{T,j}\}.$ Finally, we must assume that the set $Z^0_T$ of relative fixed points of the hamiltonian flow $\Phi^T$ on $N^*\mathcal{F}$ are clean for all $T\in \mathcal{MST}$ in the sense Definition 5. 

\begin{cor}\label{c:minstratum}
With the above assumptions,  
\begin{equation}
\chi(t)\Pi_*\Delta^*\bigl(U_B(t,x,y)\bigr)=\sum_{T\in \mathcal{MST}}\nu^0_{T}(t),
\end{equation}
where $\nu^0_{T}\in I^{-1/4-e_T/2-r}(\mathbb{R},\Gamma^T,\mathbb{R})$ where $r=-p/2.$  Furthermore, $\nu_{T}^0$ has an expansion of the form
\begin{equation}\label{e:expansionmin}
\nu_{T}^0(t)=e^{\frac{i\pi m_T}{4}}\sum_{j=0}^\infty\sigma_j(T)(t-T+i0)^{-\frac{e_T-1}{2}+\frac{p}{2}-j}\,mod\,C^\infty(\mathbb{R}),
\end{equation}
where $m_T$ is the Maslov index of $Z^0_T,$ and $\sigma_j(T)=\int_{S(Z^0_T)}d\mu_{Z_{T,j}^0}$, where for each $j\ge 0$, $d\mu_{Z_{T,j}^{0}}$ is a density on $S(Z_{T}^{0})$.
\end{cor}
Note that the regular closure result also follows from the above corollary, if one considers the $p+k_N$ dimensional foliation $(M,\overline{\mathcal{F}})$, which is non-singular under the regular closure hypothesis.
\end{section}

%
%
\begin{section}{Proof of the Main Results}

In this section we prove the spectral results presented in the previous section. We begin by analyzing the canonical relation of the Schwartz kernel of the basic projection operator. We are then in a position to prove the main theorems of Section 3.
\begin{prop}\label{p:canonical}
The canonical relation of the basic projector $P$ is given by 
\begin{equation}
\mathcal{C}=\bigcup_{k_1\le k\le k_N}\mathcal{C}_k
\end{equation}
where
\begin{equation}\label{e:cmax}
\mathcal{C}_{k_N}=\mathcal{C}_{max}=\{(\eta_y,\xi_x) \,|\,\xi_x\in N^*\overline{L},\, \overline{L}\subset \Sigma_{max},\,\eta_y\in\mathcal{L}^{max}_{\xi_x}\}
\end{equation}
and for $k<k_N$
\begin{equation}\label{e:cf}
\mathcal{C}_{k}=\{(\eta_y,\xi_x) \,|\,\xi_x\in N^*\overline{L},\, \overline{L}\subset \Sigma_{p+k},\,\eta_y\in\mathcal{L}_{\xi_x}\}.
\end{equation}
\end{prop}


\begin{proof}
First consider the canonical relation over points in the maximal stratum, $\Sigma_{k_N}=\Sigma_{max}.$  Localizing about some arbitrary $p_0\in \Sigma_{max}$, let $\chi\in C^\infty(M)$ be such that $p_0\in supp(\chi)\subset U,$ where $U\subset \Sigma_{max}$ is such that $\widehat{U}$, the saturation of $U,$ by leaf closures (which is always open) is such that $\widehat{U}$ is contained in a chain of simple distinguished open sets $U_0=U,U_1,\dots, U_\ell$ covering $\overline{L}$ where the distinguished coordinates $(x,y)$ with respect to the foliation $(\Sigma_{max},\overline{\mathcal{F}})$ are valid on each $U_i.$ Let $(x,y,\xi,\eta)$ be the corresponding coordinates on open set $V_i\subset T^*M$, with $\pi(V_i)=U_i$. If $\langle\cdot,\,\cdot\rangle_M$ denotes the distribution pairing in $M$, then
\begin{eqnarray}
\widehat{\chi P(u)}&=&\langle P(u),\, \chi \,e^{-i(x,y)\cdot(\xi,\eta)}\rangle_M\nonumber\\
&=&\langle u,\,P(\chi )\,e^{-i(x,y)\cdot(0,\eta)}\rangle_M\label{e:thing1}\\
&=&\widehat{P(\chi)u}(0,\eta), \label{e:thing2} 
\end{eqnarray}
where \eqref{e:thing1} follows from the definition of $P$ on generalized functions.  Then, the expression in \eqref{e:thing2} is rapidly decreasing if and only if there exists an open conic neighborhood $\Gamma\subset T^*M$ of $(0,\eta)$ such that $(supp(P(\chi))\times \Gamma)\cap WF(u)=\emptyset.$   

Notice that $supp(P(\chi))\subset \widehat{U}$ is saturated by the leaf closures. Thus,  if \\ $p_0=(x_0,y_0)\in supp(\chi)$ then entire leaf closure consisting of points of the form $(x,y_0)$ in local coordinates  is in $supp(P(\chi))$, since $P(\chi)$ can not distinguish between points in the same leaf closure. 
It follows that if $(x_0,y_0,0,\eta)\in WF(u)$ and $(x,y_0)\notin sing\,supp(u),$ then for any open conic neighborhood $\Gamma$ of $(0,\eta)$ $(supp(P(\chi))\times V)\cap WF(u)\not=\emptyset.$  Hence, if $(x_0,y_0,0,\eta)\in WF(u)$ then \\ $(x,y_0,0,\eta)\in WF(Pu)$ for every $x$.  (And hence, observe that $P$ is not psuedolocal, although it does not propagate the singular support beyond the saturation of the support by leaf closures.)

If we consider generalized functions with support contained in $U$, we have
\begin{eqnarray}
WF(Pu)&=&\!\!\!\{  (x_1,y_0,0,\eta_0)\,|\,(x_0,y_0,0,\eta_0)\in WF(u) \}
\nonumber\\
&=&\mathcal{C}_{max}(WF(u)), 
\end{eqnarray}
where
\begin{equation}
\mathcal{C}_{max}=\{ (x_1,y_0,0,\eta_0); (x_0,y_0,0,\eta_0)\}.
\end{equation}
The representation above must be valid in every such simple distinguished open set where the coordinates above are valid.  Following a chain of overlapping simple distinguished open sets above defined along a curve $\overline{\alpha}$ contained in the leaf closure with $\overline{\alpha}(0)=p_0=(x_0,y_0)$ and $\overline{\alpha}(1)=p_1=(x_1,y_0)$, implies that $(x_1,y_0,0,\eta_0)$ and $(x_0,y_0,0,\eta_0)$ are related by $(dh_{\overline{\alpha}}^{-1})^*$
and \eqref{e:cmax} follows. 

Consider next the canonical relation over a point in an arbitrary singular stratum, $p_0\in \Sigma_{p+k}$ for $k<k_N$. Let $\chi\in C^\infty(M)$ be such that $p_0\in supp(\chi)\subset U,$ where $U$ is a simple distinguished open set in $M$ with respect to $\mathcal{F}$, analogously as above. Let $(x,y)$ be distinguished coordinates on $U$ {\it with respect to the original foliation} $\mathcal{F}.$ Let $(x,y,\xi,\eta)$ be the corresponding coordinates on the corresponding $V_i$ open in $T^*M$, as above.  Then, for an arbitrary $u\in\mathcal{D}'(M)$ and $p=(x,y)\in U$, we have 
 \begin{equation}\label{e:thing3}
\widehat{\chi P(u)}=\langle P(u),\, \chi \,e^{-i(x,y)\cdot(\xi,\eta)}\rangle_M\nonumber\\
=\langle u,\,P(\chi )\,P(e^{-i(x,y,)\cdot(\xi,\eta)})\rangle_M. 
\end{equation}
Note that $P(e^{-i(x,y,)\cdot(\xi,\eta)})=e^{-iP((x,y)\cdot(0,\eta))}=e^{-if}$ where $f$ is of the form $z\cdot \zeta$ for some transverse variables $z$ and corresponding covectors $\zeta$ such that $f$ is a basic function. It follows that $f$ must be constant on the leaf closure through $p$.  As such, $X(f)=X(z\cdot\zeta)=0$ for all $X\in T\overline{\mathcal{F}}.$  Consequently, $\zeta\in N^*_p\overline{L}$, and 
\begin{equation}
\widehat{\chi P(u)}=\langle u,\,P(\chi )\,e^{-i z\cdot \zeta}\rangle_M=\widehat{P(\chi)u}(\zeta)
\end{equation}
is rapidly decreasing if and only if there exists an open conic neighborhood $\Gamma$ of $\zeta\in N^*\overline{L}$ such that $(supp(P(\chi))\times \Gamma)\cap WF(u)=\emptyset.$  

Suppose $\xi_{p_0}=(x_0,y_0;\zeta)\in WF(u)$ and $p_1=(x_1,y_0)$ which does not necessarily belong to $sing\,supp(u).$ Then, since $p_0$ and $p_1$ belong to the same leaf and $supp(P(\chi))$ is saturated, it follows as above that $(supp(P(\chi))\times \Gamma)\cap WF(u) \not=\emptyset$ for all open conic neighborhoods $\Gamma$ of $\zeta$. Thus if $\xi_{p_0}=(x_0,y_0,\zeta)\in WF(u)$ then $\eta_{p_1}=(x_1,y_0,\zeta)$ belongs to $WF(Pu),$ for all $x_1$.  Following a chain of simple distinguished opens sets defined along a curve $\alpha$ contained in the leaf containing $p_0$ and $p_1$, we see, as before that \eqref{e:cf} holds. 

%
\end{proof}

\begin{rem}
Note, with respect to the ``leaf diagonal" 
\begin{equation}
\mathcal{C}_{\mathcal{F}}=\{(\eta_y,\xi_x)\,|\, \eta_y\in\mathcal{L}_{\xi_x}\}=\bigcup_{\mathcal{L}_{\xi_x}\subset N^*\mathcal{F}}\,\mathcal{L}_{\xi_x}\times \mathcal{L}_{\xi_x},
\end{equation} 
we have for all $k<k_N$
\begin{equation}
\mathcal{C}_k=\bigcup_{\overline{L}\subset \Sigma_{p+k}}(\bigcup_{\mathcal{L}_{\xi_x}\subset N^*\overline{L}}  \mathcal{L}_{\xi_x}\times \mathcal{L}_{\xi_x})=\mathcal{C}_{\mathcal{F}}\cap \bigcup_{\overline{L}\subset \Sigma_{p+k}} N^*\overline{L}\times N^*\overline{L}.
\end{equation}

Note also, that if the $\Sigma_{p+k_1}$ stratum consist only of compact leaves (i. e., $k_1=0$), then $\mathcal{C}_0$ is Lagrangian.
\end{rem}

\begin{rem}
Notice that in local coordinates it is easily seen that $\mathcal{C}_{max}$ is an (immersed) Lagrangian submanifold with respect to the symplectic form $\omega_1-\omega_2$ on $ T^*M\times T^*M.$   In the vicinity of a leaf closure whose holonomy is non-trivial but finite $C_{max}$ may be immersed, rather than embedded. (Note: from \cite{Molino}, Chapter 5.4, no leaf closure in $\Sigma_{max}$ can have infinite holonomy.) Note also, that $\mathcal{C}_k$ is of dimension $2n-k$ rather than $2n$ and thus cannot be Lagrangian in $T^*M\times T^*M$.  If $k_1=0$ then $C_0$ is a (possibly immersed) Lagrangian.
\end{rem}
\begin{proof}
{\it Proof of Theorem \ref{t:sojourntimes}}. Let $K$ be the Schwartz kernel of $P$ acting on half-densities, as usual:  $K\in\mathcal{D}'(M\times M,\Omega^{1/2}_{M\times M})$ with
\begin{equation}
P(f(x_2))=\int_M K(x_1,x_2)f(x_2)=(Pf)(x_1).
\end{equation}
Observe that $P$ satisfies $P^2=P,$ by \cite{PaRi}, hence the corresponding Schwartz kernel satisfies the following relation:

\begin{equation}\label{e:idempotent}
\int_M K(x_1,x_2)\,K(x_2,y_1)=K(x_1,y_1).
\end{equation}
Thus,
\begin{equation}
\Pi_*\Delta^*\bigl(U_B(t,x_2\,y_2)\bigr)=\int_M\int_M\int_M K(x_1,x_2)\,K(x_2,y_1)U(t,x_1,y_1).
\end{equation}

It then follows from \eqref{e:idempotent} that
\begin{equation}
\Pi_*\Delta^*\bigl(U_B(t,x,y)\bigr)=\Pi_*\Delta^*\bigl(P_{x}P_{y}U(t,x,y)\bigr)=\int_M\int _MK(x,y)U(t,x,y).
\end{equation}
Note: this is exactly the situation considered in (1.6) of \cite{Z}. The Schwartz kernel $K$ of $P$ is a distribution on $M\times M$ with wave front set contained in $\mathcal{C}.$ Hence, the wave front set of $\Pi_*\Delta^*\bigl(U_B(t,x,y)\bigr)$ is estimated by the following:

\begin{eqnarray}
WF\Bigl(\Pi_*\Delta^*\bigl(U_B(t,x,y)\bigr)\Bigr)&\subset &\Lambda_t'\circ\mathcal{C}\nonumber\\
&=&\{(t,\tau)\,|\,\exists \bigl((x,\xi);(y,\eta)\bigr)\in \Lambda_t'\cap \mathcal{C} \}
\end{eqnarray}

where $\Lambda_t$ is the graph of the hamiltonian flow of the metric on $T^*M$:
\begin{equation}
\Lambda_t=\{\bigl((x,\xi);(y,\eta)\bigr)\,|\,\Phi^t(x,\xi)=(y,\eta),\,\tau=|\xi|\},
\end{equation}
and $\mathcal{C}$ is the canonical relation of the basic projector $P$.  The result is immediate. 
\end{proof}

In order to study the singularities of the wave trace further for the proof of Theorem \ref{t:regclosure}, we must investigate the clean-ness of intersection of the canonical relations of $P$ and $U(t)=e^{it\sqrt{\Delta}}.$  Accordingly, we have the following:

\begin{prop}
The components of $\mathcal{C}_{max}$ of the canonical relation of $P$ and $\Lambda$ of $e^{-it\sqrt{\Delta}}$ intersect cleanly.  If $k_1=0$, then $\mathcal{C}_0$ intersects $\Lambda$ cleanly also.
\end{prop}
\begin{proof}
We will show the second part of the proposition first. Let $Y$ denote the set
$T^*M\times \Delta(T^*M\times T^*M)\times T^*M\times T^*\mathbb{R}$ where $ \Delta(T^*M\times T^*M)$ denotes the diagonal in the product space. 

In the notation of Hormander's clean intersection criteria, \cite{H3}, we consider the set
\begin{equation*}
(\mathcal{C}_0\times \Lambda' )\cap Y.
\end{equation*}
Observe first that as a consequence of the discussion of Section 2.2, the flow restricts to $N^*\mathcal{F}.$  Hence, the canonical relation of $e^{-it\sqrt{\Delta}},$ which is just the graph of the hamiltonian flow restricts to $N^*\mathcal{F}\times N^*\mathcal{F}\times T^*\mathbb{R}.$  Let $\Lambda_ {N^*\mathcal{F}}$ denote $\Lambda\cap N^*\mathcal{F}\times N^*\mathcal{F}\times T^*\mathbb{R}.$ If $pr_1$ and $pr_2$ denote the projections from $T^*M\times T^*M\rightarrow T^*M$ onto the first and second components, then observe that $pr_2(\mathcal{C}_0)=N^*\mathcal{F}.$  Thus,
\begin{equation}\label{e:intersection0}
(\mathcal{C}_0\times \Lambda' )\cap Y=(\mathcal{C}_0\times \Lambda' _{N^*\mathcal{F}})\cap Y.
\end{equation}
Let $p\in (\mathcal{C}_0\times \Lambda' _{N^*\mathcal{F}})\cap Y$.  Then $p$ has the form 
\begin{equation}\label{e:typicalp}
p=\bigl((xi_x;\Phi^t(\eta_y);\Phi^t(\eta_y);\eta_y;(t,\,\tau)\bigr),
\end{equation}
where $(\xi_x;\Phi^t(\eta_y))\in \mathcal{C}_0\subset N^*\mathcal{F}\times N^*\mathcal{F}$ with $\xi_x=(dh_{\gamma}^{-1})^*\Phi^t(\eta_y)$ for some suitable $\boldsymbol{\gamma}\in \mathcal{G}(\mathcal{F})$ and $\tau>0$. The set of such points is a manifold because it is equal to the set $\mathcal{C}_0\times N^*\mathcal{F} \times S$ where $S$ is the inverse image of $N^*\mathcal{F}$ by the restricted flow, which is a diffeomorphism. 

Now consider the tangent space to the intersection $\mathcal{C}_0\times \Lambda'_{N^*\mathcal{F}}\cap Y.$ For the intersection to be clean, then for all points $p$ belonging to the intersection
\begin{equation}
T_p(\mathcal{C}_{0}\times \Lambda'_{N^*\mathcal{F}}) \cap T_pY\subset T_p((\mathcal{C}_{0}\times \Lambda'_{N^*\mathcal{F}})\cap Y).
\end{equation}
To verify this, we first characterize the tangent space of the intersection.  With respect to the splitting given by the null foliation on $N^*\mathcal{F}$, any tangent vector $X\in T(N^*\mathcal{F})$ splits into $X'+X''$ where $X'\in T\widetilde{\mathcal{F}}$ and $X''=H(X)\in N\widetilde{\mathcal{F}}.$
Consider $T_{p_1}\mathcal{C}_{0},$ where $p_1=(\xi_x,\Phi^t(\eta_y))$ where $p$ in \eqref{e:typicalp} equals $(p_1,p_2).$ Then, from \eqref{e:liftedhol}, if $((X',X''), (Z',Z''))\in T_{p_1}\mathcal{C}_{0},$
then
\begin{equation}
X''=d\widetilde{h}_{(\gamma,\Phi^t(\eta_y))}(Z'')
\end{equation}
for some suitable $\gamma\in \mathcal{G}(\mathcal{F}).$  Recalling Lemma \ref{l:phlow}, $d\Phi^t$ splits with respect to the splitting of the tangent space of $N^*\mathcal{F}$ and we see that at a point $p$ belonging to the intersection \eqref{e:intersection0}, the tangent space consists of vectors of the form:
\begin{equation}\label{e:clean0}
\bigl( (X',d\widetilde{h}_{(\gamma,\Phi^t(\eta_y))}(d\Phi^t(Y'')),(d\Phi^t(Y'),d\Phi^t(Y''), (d\Phi^t(Y'),d\Phi^t(Y''),Y',Y'', W\bigr)
\end{equation}
where $\bigl((d\Phi^t(Y'),d\Phi^t(Y''),(Y',Y''), W\bigr)\in T_{p_2}\Lambda'_{N^*\mathcal{F}}.$

Now consider a vector in $T_p(\mathcal{C}_{0}\times \Lambda'_{N^*\mathcal{F}}) \cap T_pY.$ Such a vector must be of the form $\bigl((X',X''),(Y',Y''),(Y',Y''),(Z', Z''), W)\bigr)$ where $\bigl((X',X''),(Y',Y'')\bigr)\in T_{p_1}\mathcal{C}_0$ and $\bigl((Y',Y''),(Z', Z''), W)\bigr)\in T_{p_2}\Lambda'_{N^*\mathcal{F}}.$  But then $Y'=d\Phi^t(Z')$ and $Y''=d\Phi^t(Z'')$ and, hence  for some $\gamma\in\mathcal{G}(\mathcal{F}),$ $X''=d\widetilde{h}_{(\gamma,\Phi^t(\eta_y))}(d\Phi^t(Z''))$, and we have a vector of the form \eqref{e:clean0}, proving clean-ness.

For the first part of the proposition, the reasoning is similar, with some modifications. As before, observe that $pr_2(\mathcal{C}_{max})=N^*\overline{\mathcal{F}}_{max}\subset N^*\mathcal{F}$ is a submanifold, and
\begin{equation}\label{e:intersection}
(\mathcal{C}_{max}\times \Lambda' )\cap Y=(\mathcal{C}_{max}\times \Lambda' _{N^*\mathcal{F}})\cap Y.
\end{equation}
Let $p\in (\mathcal{C}_{max}\times \Lambda' _{N^*\mathcal{F}})\cap Y$.  Then $p$ has the form $$p=\bigl(\xi_x;\Phi^t(\eta_y);\Phi^t(\eta_y);\eta_y;(t,\,\tau)\bigr),$$ 
where $(\xi_x;\Phi^t(\eta_y))\in \mathcal{C}_{max}\subset N^*\overline{\mathcal{F}}\times N^*\overline{\mathcal{F}}$ with $\xi_x=(dh_{\overline{\gamma}}^{-1})^*\Phi^t(\eta_y)$ for some suitable $\boldsymbol{\overline{\gamma}}\in \mathcal{G}_{max}$ and $\tau>0$.

The set of such points forms a manifold because it is precisely equal to the set
$\mathcal{C}_{max}\times pr_2(\mathcal{C}_{max})\times S$ where $S$ is the inverse image via a diffeomorphism of $N^*\overline{\mathcal{F}}_{max}$ by the flow restricted to $N^*\mathcal{F}.$  Observe that the intersection is a product of manifolds.

Next, consider $T_{p_1}\mathcal{C}_{max},$ where $p_1=(\xi_x,\Phi^t(\eta_y)).$ Recall that $N^*\overline{\mathcal{F}}_{max}$ is a saturated manifold by $T\widetilde{\mathcal{F}}_{max}$ and, as in Section 2.2,\\ $T(N^*\overline{\mathcal{F}}_{max})=T\widetilde{\mathcal{F}}_{max}\oplus N\widetilde{\mathcal{F}}_{max}$ with $H_{max}$ the projection onto the horizontal space.  Then the holonomy relation $(x,\xi)=(dh_{\overline{\gamma}}^{-1})^*\Phi^t(y,\eta)$ for suitable $\boldsymbol{\overline{\gamma}}\in \mathcal{G}_{max}$ implies that if $(X,d\Phi^t_{\eta_y}(Y))\in T_{p_1}\mathcal{C}_{max},$ then
\begin{equation}\label{e:barholonomy}
H_{max}(X)=d\widetilde{h}^{max}_{(\overline{\gamma},\Phi^t(\eta_y))}(H_{max}(d\Phi^t_{\eta_y}(Y)))
\end{equation}

Thus, at a point $p$ belonging to the intersection \eqref{e:intersection}, the tangent space consists of vectors of the form:
\begin{equation}\label{e:clean}
\bigl( X,d\Phi^t(Y),d\Phi^t(Y),(Y,Z)\bigr)
\end{equation}
where $(X,d\Phi^t(Y))$ satisfies \eqref{e:barholonomy}.  The rest of the clean-ness argument goes through, as in the first case considered above.

\end{proof}

\begin{proof} 
{\it Proofs of Theorem \ref{t:partialtrace} and Corollary \ref{t:regclosure}.} Note that the corollary follows immediately if $\mathcal{SST}=\emptyset.$  However, the direct proof of the corollary is also instructive. If the foliation has regular closure, then there is only one component of the canonical relation, and thus $\mathcal{C}$ is a manifold, by the previous proposition.  In this case, the composition of the entire canonical relation of $P$ and $U$ is clean, and we can compute the trace of the basic wave kernel at non-zero relative periods $T$ by the standard stationary phase arguments and clean intersection arguments as in, for example, the proof of Proposition 1.10 of \cite{Z}, or the corresponding result of \cite{DG}.  

Let $\Gamma$ be the canonical relation of the Schwart kernel of $\Pi_*\Delta^*,$ the conormal to the ``diagonal" in $T^*(\mathbb{R}\times M\times M\times \mathbb{R}),$ let $\Lambda^B=\mathcal{C}\circ \Lambda'_{N^*\mathcal{F}}$, and let $F$ denote the fibre product $\{(s,t)\in \Gamma'\times \Lambda^B\,|\,f(s)=\iota(t)\}.$  Let $f$ denote the embedding $f:\Gamma'\rightarrow T^*(\mathbb{R}\times M\times M)$ with $f(t,\tau,x,\xi,x,-\xi,t,-\tau)=(t,\tau,x,\xi)$.  The clean intersection theory of \cite{DG} may be applied if the fibre product diagram below is clean:
\begin{equation}
\xymatrix{ \Gamma'\ar[d]_f& F\ar[l]_{p_1}\ar[d]^{p_2}\\
T^*(\mathbb{R}\times M\times M) &\ar[l]_\iota\mathcal{C}\circ \Lambda'_{N^*\mathcal{F}}}
\end{equation}
Note that $Z:=\cup_{T\in \mathcal{ST}}S(Z^{max}_T)$ is the compact fibre of the map \\ $\xymatrix{F\ar[r]^{p_1}& \Gamma'\ar[r] &\Gamma'\circ \Lambda^B},$ where the second arrow is projection onto the last component of $\Gamma'$, $T^*\mathbb{R}.$  The diagram will be clean if each $Z^{max}_T$ is a manifold and the associated diagram below is also a fibre product
\begin{equation}
\xymatrix{ T_x\Gamma'\ar[d]_{df_x}&T_p F\ar[l]_{dp_1}\ar[d]^{dp_2}\\
T_z^*(\mathbb{R}\times M\times M) &\ar[l]_{d\iota_y} T_y(\Lambda^B)'
}
\end{equation}
where $p=(x,y)$ and $z=f(x)=\iota(y)$.
By splitting up the tangent space to $Z^{max}_T$ into horizontal and leafwise parts, we see that this occurs precisely when the clean-ness condition of Definition 4 is satisfied for $Z^{max}_T$. 

If the set of relative fixed points $Z^{max}_T$, is clean, then it is possible to define smooth positive densities on $T(Z^{max}_T)$, denoted by $d\mu_{Z^{max}_{T,j}},$ as follows.  First, recall the characteristic form of a foliation of dimension $p$, $\chi_{\mathcal{F}}:$ let $X_1, \dots, X_p$ be in $T\mathcal{F}$, and let $\{E_j\}_{j=1}^n$ be an orthonormal frame of $M$, of dimension $n$.  Then define the canonical $p$-form via the metric on $M$:
\begin{equation}
\chi_{\mathcal{F}}(X_1,\dots,X_p)=det(g_{ij}(E_i,X_j)).
\end{equation}
Applying this to the $p+k_N$ dimensional foliation $(\Sigma_{max},\overline{\mathcal{F}})$, we can define a canonical leafwise density on $T\widetilde{\mathcal{F}}^{max}$ by lifting the to $(T\mathcal{F}_S)^0,$ via $\pi:T^*M\rightarrow M$.  Thus, $|d\pi^*\chi_{\overline{\mathcal{F}}}|$ defines a positive leafwise density on  $(T\overline{\mathcal{F}})^0.$  Next, since the horizontal space, $N\widetilde{\mathcal{F}}^{max}$ of $((T\overline{\mathcal{F}})^0,\,T\widetilde{\mathcal{F}}^{max})$ is a symplectic space, and $\Phi_T$ and $d\widetilde{h}_{\boldsymbol{\overline{\alpha}}}$ are symplectic diffeomorphisms of  $N\widetilde{\mathcal{F}}^{max}$, one can use Section 4 of \cite{DG} to construct a canonical densities on $N\widetilde{\mathcal{F}}$, say $d\mu'_{Z_{T,j}^{max}}$.  One then constructs the densities on each component of $Z_T^{max}$ as follows:
\begin{equation}
d\mu_{Z_{T,j}^{max}}=d\pi^*(\chi_{\overline{\mathcal{F}}})\otimes d\mu'_{Z_{T,j}^{max}}.
\end{equation}
Once we have the canonical densities on $Z_T^{max}$, we obtain densities on $S(Z_T^{max})$ in the usual way, (see Section 4, \cite{DG}).

To compute the order $r$, observe once more that $P^2=P.$ Furthermore, its kernel $K$ belongs to $I^r(M\times M,\mathcal{C}),$ hence, $P^2\in I^r(M\times M,\mathcal{C}).$  By H\"{o}rmander's composition theorem for such distributions , $P^2\in I^{2r+e/2}(M\times M,\mathcal{C}'\circ\mathcal{C}')$ where $e$ is the excess in the composition (which is clean) $\mathcal{C}'\circ\mathcal{C}'.$  Thus $r=-e/2.$ The excess $e$ is the dimension of the fibre of the projection 
\begin{equation}
\mathcal{C}\times \mathcal{C}'\cap (T^*M\times \Delta(T^*M\times T^*M)\times T^*M)
\end{equation}
which is $p+k_N.$

The calculation of the leading order part arises from the calculation of the symbol of $K$, which results in 
\begin{equation}
\sigma_K=pr^*|dx\wedge dy \wedge d\zeta|^{1/2}
\end{equation}
where $pr:N^*(\Delta(M\times M))\rightarrow N^*\overline{\mathcal{F}}$ and $(x,y,\xi,\eta)$ are coordinates on $T^*M$ defined by distinguished coordinates with respect to the foliation $(\Sigma_{max},\overline{\mathcal{F}})$.  (This is just pull-back of the volume half-density on the conormal bundle of the leaf closure.)  These coordinates can be defined on all of $T^*M$ since the leaf closures have constant dimension under the hypotheses of the corollary.  

To prove Theorem \ref{t:partialtrace}, observe that with the clean-ness condition of Definition 4, we can still apply clean intersection theory to the cut-off wave trace $\chi(t)\Pi_*\Delta^*PU(t,x,y)$ since this only involves the composition of the $\mathcal{C}_{max}$ component of the canonical relation of $P$ with $\Lambda$ for $T\in supp(\chi).$  In other words, if we let $\Lambda_T=\{\bigl((x,\xi);(y,\eta)\bigr)\,|\,\Phi^T(x,\xi)=(y,\eta),\,\tau=|\xi|\}$, where $\Phi^t$ is the restricted flow, and let $\Lambda_T^{max}=\mathcal{C}_{max}\circ \Lambda_T',$ then the we only require clean-ness for $\mathcal{C}_{max}$ as follows: 
\begin{equation}
\xymatrix{ \Gamma'\ar[d]_f& F\ar[l]_{p_1}\ar[d]^{p_2}\\
T^*(\mathbb{R}\times M\times M) &\ar[l]_\iota(\Lambda_T^{max})'}
\end{equation}
As before, the fibre of $\xymatrix{F\ar[r]^{p_1}& \Gamma'\ar[r] &\Gamma'\circ \Lambda^B},$ is $S(Z_T^{max})$, and the diagram will be clean if $Z_T^{max}$ is a manifold and the associated diagram below is also a fibre product
\begin{equation}
\xymatrix{ T_x\Gamma'\ar[d]_{df_x}&T_p F\ar[l]_{dp_1}\ar[d]^{dp_2}\\
T_z^*(\mathbb{R}\times M\times M) &\ar[l]_{d\iota_y} T_y(\Lambda_T^{max})'
}
\end{equation}
which is satisified when the clean-ness condition of Definition 4 holds.

Since $\Sigma_{max}$ is an open dense set in $M$, one can locally represent the composition of Schwartz kernel of $P$ and $U(t,x,y)$ as a locally finite sum of integrals in local coordinates over $\Sigma_{max}.$ The techniques used in the analysis of the cut-off wave trace are the usual stationary phase arguments applied to these integral expressions.  These arguments are entirely local in $t,$ and one can simply perform the usual stationary phase arguments locally in $t$ on the component of the canonical relation of $P$ that corresponds to the fixed points corresponding to $T$ which, by hypothesis, are associated only to the maximal stratum.
\end{proof}

Finally, the proof of Corollary \ref{c:minstratum} follows in an analogous manner to the proof of Theorem \ref{t:partialtrace}, using Definition 5 in place of Definition 4.
\end{section}

%
%
\begin{section}{Examples}
In this section we present three examples that illustrate some of the behavior of the transverse geometry. The first example has regular closure, and the second and third do not. All of the examples are non-simple foliations generated by suspensions, which are a bit special in the class of foliations.  They have a basic mean curvature equal to zero, a global transversal manifold, and proper leaves, which are totally geodesic.  Furthermore, the metric is a product metric, and thus the functions $H_{\mathcal{F}}$ and $H_{\mathcal{F}^\perp}$ Poisson commute. Relatively closed curves, in this case, are quite likely closed in the ordinary sense.  For such examples, it seems likely that a better wave trace formula may be possible.   
\begin{subsection}{The Suspension of an Irrational Rotation of a Torus}
Consider  the suspension of an irrational rotation about the $z$-axis on the 2-torus with the usual round metric.  Let $(\psi,\theta,s)$ be coordinates on $T^2\times [0,1],$ with $0\le\psi<2\pi,$ $0\le\theta<2\pi,$ $0\le s\le 1.$  Let $\alpha$ be irrational multiple of $2\pi$, and define a $\mathbb{Z}=\pi_1(S^1)$ action on $T^2\times [0,1]$ by rotation by $k\alpha$ about the $z$-axis in the $T^2$ component and translation by $k$ on the $[0,1]$ component:  
$k\cdot (\psi,\theta,0)=(\psi, \theta+k\alpha,s+k (\text{mod } 1)).$  Note that this action has no fixed points.
Our foliated manifold is $M$ is $T^2\times [0,1] / \sim $ where $(\psi,\theta,0)\sim (\psi,\theta+\alpha,1).$

The metric on this manifold will be the usual product metric on $T^2\times [0,1]$.  The coordinates above are orthogonal, and the facts that (1) the rotation in the $\theta$ coordinate is an infinitesimal isometry, and (2) that the leaves are totally geodesic will imply that the metric depends only on $\psi,$ and the mean curvature form is zero. 

The leaves of this foliation are the one dimensional submanifolds:
\begin{equation}
L_{(\psi,\theta)}=\bigcup_{k\in \mathbb{Z}}\{\psi\}\times \{\theta+k\alpha\}\times [0,1]
\end{equation}
and the leaf closures are the two dimensional submanifolds:
\begin{equation}
\overline{L_{\psi}}=\{(\psi,\theta,s)\,|\,0\le\theta<2\pi,\,s\in[0,1]\}.
\end{equation}

Thus, none of the leaves will be closed, and so the foliation is not simple, but the leaf closures all have the same dimension. Hence, this example will have regular closure, and the partition of the manifold into leaf closures will be another  foliation of the manifold.  In fact, the foliation by leaf closures of the original foliation is a simple foliation. 

Basic functions for this example consist of the functions of $\psi:$ if one considers a cube in $(\psi,\theta,s)$ coordinates of the form $[0,2\pi]\times [0,\,2\pi)\times [0,1]$ where $(\psi,\theta,0)\sim (\psi,\theta+\alpha,1)$, one sees that the only continuous functions that are constant on the leaves must also be constant in $\theta$. Thus, the smooth basic functions for this foliation are just the functions of the $\psi$ variable that are smooth on $S^1$, due to the identification of the points $(z,\theta, 0)$ and $(z,\theta+\alpha,1)$. 

The basic Laplacian for this example is just the Laplacian on $S^1$, induced from the round metric on $T^2$ on the longitudinal circle. Hence, the basic spectrum of the Laplacian is, of course, just $\{k^2\}$ where $k\in \mathbb{N}.$

If one considers the sojourn times for this example, they are a discrete set corresponding to the multiples of the lengths of the meridian circles on the torus.  

Note that in this example, all of the leaf closures have trivial holonomy, so here the space of leaf closures is quite tame--it is actually a manifold, rather than an orbifold.
\end{subsection}

\begin{subsection}{The Suspension of an Irrational Rotation of a Sphere}
Now consider an analogous example to the one above by repeating the construction with the sphere in place of the torus. Endow the 2-sphere $S^2$ with the usual round metric, and cylindrical coordinates $(z,\,\theta)$ for  $-1\le z\le 1,\,0\le\theta<2\pi$.  Now let $\alpha$ be an irrational multiple of $2\pi.$  Let $(z,\,\theta,\, s)$ be coordinates on $S^2\times [0,1].$ Our manifold $M$ will be $S^2\times [0,1]/ \sim$ where $(z,\,\theta,\,0)\sim(z,\,\theta+\alpha,\,1).$ $M$ is the orbit space of a $\mathbb{Z}=\pi_1(S^1)$ action on $S^2\times [0,1]$ which is defined by a rotation by $k\alpha$ on the $S^2$ component and by translation by $k$ on the last component.  Observe that this action has fixed points at the poles of $S^2$. (Note:  This example appears in Section 4 of \cite{KR2}.  There it is shown explicitly for this example that the heat kernel, $K(t,x,x),$ is not integrable over $M$.)

The leaves of this foliation, $\mathcal{F},$ are the one-dimensional submanifolds, indexed by the points on $S^2$:
\begin{equation}
L_{(z\,,\theta)}=\bigcup_{k\in\mathbb{Z}} \{ z \}\times \{\theta +k\alpha\}\times [0,1].
\end{equation}

We observe that this foliation is not simple because the leaves are not the connected components of the inverse images of a smooth submersion on $M$.  In particular, the leaves are not closed, except for the leaves in $M$ that over the North and South poles in the $S^2$ component ($z=\pm 1$).  

Observe that here there are two types of leaf closures:

\begin{eqnarray}
\overline{L_z}&=&\{(z,\,\theta,\,s) |\, 0\le\theta<2\pi,\,0\le s\le 1\}, \text { for } |z|<1 \label{fatclosedleaf}\\
\overline{L_x}&=&\{(x,\,s) |\, 0\le s\le 1\},\,x\in\{N,S\}\label{thinclosedleaf}
\end{eqnarray}
where $N$ and $S$ refer to the north and south poles on $S^2$.
Notice that the leaf closures in \eqref{fatclosedleaf} are of dimension 2, while the leaf closures over the poles in \eqref{thinclosedleaf} are of dimension 1.  

In terms of the discussion of the partition of the manifold into strata, we see that $k_N=0$ and $k_N=1$ and thus there are two strata:
\begin{eqnarray}
\Sigma_1&=&\overline{L_N}\cup \overline{L_S},\\
\Sigma_2&=&\bigcup_{|z|<1}\overline{L_z}.
\end{eqnarray}

The metric on $M$ will be the usual product metric so that $(z,\,\theta,\,s)$ is an orthogonal coordinate system, and the metric is given by 

\begin{eqnarray}\label{e:metric}
g_{11}&=&\langle \partial_z,\,\partial_z\rangle=\frac{1}{(1-z^2)}\nonumber \\
g_{22}&=&\langle \partial_\theta,\,\partial_\theta\rangle=1-z^2\nonumber\\
g_{33}&=&\langle \partial_s,\,\partial_s\rangle=1,\nonumber\\
g_{ij}&=&0 \text{ for }i\not= j.
\end{eqnarray}

If one considers a cube in $(z,\theta,s)$ coordinates of the form $(-1,\, 1)\times [0,\,2\pi)\times [0,1]$ where $(z,\theta,0)\sim (z,\theta+\alpha,1)$, one sees that the only continuous functions that are constant on the leaves must also be constant in $\theta$. Thus, the smooth basic functions for this foliation are just the functions of the $z$ variable that are smooth on $(-1,\,1)$ and continuous on the closure of this interval, due to the identification of the points $(z,\theta, 0)$ and $(z,\theta+\alpha,1)$.  

Observe from \eqref{e:metric} that the Christoffel symbols for this metric depend only on $z.$ Notice that this implies that the mean curvature form $\kappa$ is basic.    (The vanishing of $\kappa$ is related to the fact that this foliation is, in fact, totally geodesic in this metric.)  So we have a non-simple example of a foliation whose mean curvature is basic. 

Note that in this metric, the Laplacian and the basic Laplacian on $C^\infty(M)$ have the following expressions in $(z,\theta,s)$ coordinates:

\begin{eqnarray}
\Delta&=&-(1-z^2)\partial^2_z+2z\partial_z-\frac{1}{1-z^2}\partial^2_\theta-\partial_s^2,\\
\Delta_B&=&-(1-z^2)\partial^2_z+2z\partial_z.
\end{eqnarray}

With respect to the basic Laplacian, it can be shown (see for example Section 4 of \cite{KR2}) that the basic spectrum is just the spectrum of the Laplacian on $S^2:$  $\{k(k+1)\},\, k\in\mathbb{N}$ with multiplicity $(2k+1).$ 

To understand the different nature of the case of a foliation with leaf closures of variable dimension, consider the lifted foliation on $\widehat{M}.$ Recall, that the basic projector $P$ is defined in terms of the operator $A$ that averages over these leaf closures. In this example, $\widehat{M}$ is an $SO(2)$ bundle.  Away from the poles, $\hat{x}$ can be denoted by the coordinates $(z,\theta, s,\phi), $ where $\phi$ denotes the coordinate on $SO(2)\cong S^1$.  There are two types of leaf closures for the lifted Foliation:
\begin{eqnarray}
\overline{K_{(z,\phi)}}&=&\{(z,\,\theta,\,s,\,\phi) |\, 0\le\theta<2\pi,\,0\le s\le 1\},\nonumber \\
&\quad&\text { for } |z|<1, \, 0\le \phi < 2\pi\label{fatclosedleafinmhat}\\
\overline{K_{x}}&=&\{(x,\,s,\,\phi) |\, 0\le s\le 1, 0\le \phi < 2\pi\}, \,x\in\{N,S\}.\label{thinclosedleafinmhat}
\end{eqnarray}
Each of these leaf closures in $\widehat{M}$ has the structure of a principal subbundle over the corresponding leaf closure in $M$, but the structure group varies depending on the leaf closure.  For the leaf closures in \eqref{fatclosedleafinmhat} the structure group is $H=\{e\},$ while in \eqref{thinclosedleafinmhat}, the structure group is $H=SO(2).$  


In our applications, we are interested in the sojourn times for this foliation. In this example, the sojourn times $\mathcal{ST}=\mathcal{MST}$ and $\mathcal{RST}=\emptyset$, and one may apply Corollary \ref{c:minstratum}. The relatively closed hamiltonnian curves correspond to multiples of $2\pi,$ the length of the meridian circles on the sphere.  
\end{subsection}

\begin{subsection}{The Suspension of an Irrational Rotation of the Cartesian Product of an Arbitrary Manifold and Sphere}

Let $X$ be any compact manifold.  We by repeat the construction of the previous example with $X\times S^2$ in place of $S^2.$ Endow $X$ with any metric and the 2-sphere $S^2$ with the usual round metric, and cylindrical coordinates $(x,z,\,\theta)$ where $x$ are local coordinates on $X$, and $-1\le z\le 1,\,0\le\theta<2\pi$.  Let $\alpha$ be an irrational multiple of $2\pi$ as before and let $(x,z,\,\theta,\, s)$ be coordinates on $(X\times S^2)\times [0,1].$ Our manifold $M$ will be $(X\times S^2)\times [0,1]/ \sim$ where $(x,z,\,\theta,\,0)\sim(x,z,\,\theta+\alpha,\,1).$  

This example is essentially the same as the previous one, except the codimension is greater.  In particular, there are now regular sojourn times corresponding to the lengths of closed geodesics that remain inside the regular stratum.  This includes closed geodesics in $X$.  If $X$ is such that the length of these relatively closed geodesics can be separated appropriately from the singular sojourn times of the previous example, then one can apply both Theorem \ref{t:partialtrace} and Corollary \ref{c:minstratum}.  However, note that the singularities arising from Corollary \ref{c:minstratum} will be of higher order.
\end{subsection}

\end{section}

%

\end{document}